\documentclass[a4paper]{amsart}

\usepackage{amsmath}   
\usepackage{amsthm}   

\usepackage{amssymb}
\usepackage[PostScript]{diagrams}
\diagramstyle[labelstyle=\scriptstyle]
\def\CD{\diagram[amstex]}

 \theoremstyle{plain} 
 
       \newtheorem{theorem}{Theorem}[section]
       \newtheorem{proposition}[theorem]{Proposition} 
       \newtheorem{corollary}[theorem]{Corollary}

 \theoremstyle{definition} 
       \newtheorem{remark}[theorem]{Remark}  
       \newtheorem{definition}[theorem]{Definition}  
       \newtheorem{example}[theorem]{Example}     

\newcommand{\GGAmi}{A_{\min}}
\newcommand{\GGAma}{A_{\max}}
\newcommand{\GGinj}{\text{\rm i}}
\newcommand{\GGwA}{\widetilde A}
\newcommand{\GGdd}{\text{\it \dj\hspace{1pt}}}
\newcommand{\GGOp}{\operatorname{Op}}
\newcommand{\GGpr}{\operatorname{pr}}
\newcommand{\GGcomega}{\overline\Omega }
\newcommand{\GGsimto}{\overset\sim\rightarrow}

\newcommand{\GGrn}{{\mathbb R}^n}
\newcommand{\GGrnp}{{\mathbb R}^n_+}
\newcommand{\GGrp}{ \mathbb R_+}

\begin{document}

\title
[Extension theory with
pseudodifferential methods]
{Extension theory for\\ elliptic partial differential operators \\with pseudodifferential methods}

\author {Gerd Grubb} 
\address
{Department of Mathematical Sciences,
Universitetsparken 5, DK-2100 Copenhagen, Denmark,
E-mail {\tt grubb@math.ku.dk}}
\begin{abstract}
This is a short survey on the connection between general extension
theories and the study of realizations of elliptic operators $A$ 
on smooth domains in $\mathbb R^n$, $n\ge 2$. The theory of 
pseudodifferential boundary problems has turned out to
be very useful here, not only as a formulational framework, but also for the
solution of specific questions. We recall some elements of that
theory, and show its application in several cases (including new
results), namely to the lower boundedness question, and the question
of spectral asymptotics for differences between resolvents.
\end{abstract}

\maketitle

\section {Introduction}\label{GG1intro}

The general theory of
 extensions of a symmetric operator (or a dual pair of
operators) in a Hilbert space, originating in the mid-1900's, has
been applied in numerous works to ordinary differential equations
(ODE), and also in a (smaller) number of works to partial
differential equations (PDE). 

There is a marked difference between the
two cases: In ODE, the playground for boundary conditions is usually
finite-dimensional vector spaces, where linear conditions can be
expressed by the help  of matrices. Moreover, the domains of
differential 
operators
defined by closure in $L_2$-based Hilbert spaces  
can usually all be
expressed in terms of functions with the relevant number of absolutely
continuous derivatives. 

In contrast, boundary conditions for PDE (in space dimensions
$n\ge 2$) are prescribed on infinite-dimensional vector spaces.
Moreover, the domains of differential operators in $L_2$-based spaces
will contain functions
with distribution derivatives, not continuous and possibly highly irregular.

Whereas extensions of ODEs can usually be described in terms of
matrices, the tools to interpret extensions in PDE cases are therefore
much more complicated. We shall give a survey of some tools developed
through the years, and their applications, emphasizing the use of
pseudodifferential operators.

\medskip

\noindent{\it Outline.} In Section \ref{GG2bvp} we recall the basic issues of
elliptic boundary value problems. Pseudodifferential operators ($\psi
$do's) are
introduced in Section \ref{GG3psdo}, and Section \ref{GG4psdbo} introduces pseudodifferential
boundary operators ($\psi $dbo's). In Section \ref{GG5ext} we recall the elements of a general abstract
extension theory, and in Section \ref{GG6impl} we show how this is
implemented for 
realizations $\GGwA$ of an elliptic operator $A$ on a domain 
$\Omega\subset \mathbb R^n$. Section \ref{GG7res} focuses on the resolvent 
formulas that can be obtained via the general theory. In
the last sections we go through several cases where pseudodifferential
techniques have proved extremely useful (some of the results here are
quite recent): In Section \ref{GG8appl1} it is the
question of whether lower boundedness holds simultaneously for a
realization $\GGwA$ and the operator $L$ over the boundary that enters in
the corresponding boundary condition; the new results deal
with unbounded domains. In Section \ref{GG9appl2} it is the
question of showing Weyl-type spectral asymptotics formulas for
differences between resolvents. The results there go back to the early
theory, and Section \ref{GG10new} presents some
additional new results.

\section {Elliptic boundary value problems} \label{GG2bvp}

In the following we use the customary multi-index notation for
differential operators: $\partial=\partial_x=(\partial_1,\dots,\partial_n)$,
$\partial_j=\partial_{x_j}=\partial/\partial x_j$, and 
$D=D_x=(D_1,\dots,D_n)$, $D_j=D_{x_j}=-i\partial/\partial x_j$; then
$\partial^\alpha =\partial_1^{\alpha _1}\cdots \partial_n^{\alpha
_n}$, $D^\alpha =D_1^{\alpha _1}\cdots D_n^{\alpha _n}$, for $\alpha
\in {\mathbb N}_0^n$; here $|\alpha |=\alpha _1+\dots+\alpha _n$.

A differential operator of order $m>0$,
\begin{equation*}
A=\sum_{|\alpha |\le m}a_\alpha (x)D^\alpha 
\end{equation*}
is said to be {\it elliptic}, resp.\ {\it strongly elliptic}, on a set
$U\subset{\mathbb R}^n$, when the principal symbol 
\begin{equation*}
a_m(x,\xi )= \sum_{|\alpha |= m}a_\alpha (x)\xi ^\alpha 
\end{equation*}
 satisfies
\begin{equation*}
a_m(x, \xi )\ne 0,\text{ resp.\ }\operatorname{Re}a_m(x, \xi )> 0,
\end{equation*}
for $x\in U$, $\xi \in{\mathbb R}^n\setminus\{0\}$.
As a basic example, the Laplacian $\Delta =\partial_1^2+\dots+\partial_n^2$
has principal symbol (and symbol) equal to $-|\xi |^2$, so it is
elliptic, and $-\Delta $ is strongly elliptic. The Laplacian has been
studied for several hundred years, and the problems around it solved
by explicit solution formulas. It is the cases with variable
($x$-dependent) coefficients, and domains more general than simple
geometric figures,  that have been a challenge in more
modern times.

The problem
\begin{equation}
Au=f\label{GG2.4}
\end{equation}
for a given function $f$ on a subset $\Omega $ of ${\mathbb R}^n$ usually has
infinitely many solutions.  To get a problem with unique
solvability, we must adjoin extra conditions such as suitable boundary
conditions. We can consider $A$ with a domain consisting of
the functions satisfying the boundary condition, as an operator $\GGwA$ acting
between suitable spaces.

Then the question of {\it existence} of a solution corresponds to the question of
whether $\GGwA$ is {\it surjective}, and the question of {\it uniqueness} of a
solution corresponds to the question of whether $\GGwA$ is
{\it injective}. In this way, the question of solvability of differential
equations is turned into a question of properties of specific
operators. The operator point of view became particularly fruitful when it was combined with
appropriate scales of function spaces, such as the {\it Sobolev spaces}, 
Sobolev \cite{GGSo50},
and with 
{\it Distribution theory}, Schwartz \cite{GGSch50}. 

When $\Omega $ is a smooth open subset of ${\mathbb R}^n$ with boundary
$\partial\Omega =\Sigma $, we refer to the standard
$L_2$-Sobolev spaces, with the following notation: $H^s({\mathbb R}^n)$
($s\in 
{\mathbb R}$) has the norm
$\|v\|_s
=\|{\mathcal F}^{-1}(\langle\xi \rangle^s{\mathcal
F}v)\|_{L_2({\mathbb R}^n)}$; here ${\mathcal F}$ is the Fourier transform 
\begin{equation*}
\mathcal F\colon u(x)\mapsto (\mathcal Fu)(\xi )=\hat u (\xi )=\int _{\GGrn}e^{-i x\cdot\xi }u(x)\, dx,
\end{equation*}
and
$\langle\xi \rangle=(1+|\xi |^2)^{\frac12}$. Next, with $r_\Omega $ denoting restriction to $\Omega $, 
\begin{equation*}
H^s(\Omega )=r_{\Omega }H^s({\mathbb
R}^n),\end{equation*} 
provided with
the norm $\|u\|_{s}=\inf\{\|v\|_s\mid v\in H^s({\mathbb
R}^n),\, u=r_\Omega v \}$. Moreover,   
\begin{equation*}
H^s_0(\GGcomega )
=\{u\in H^s({\mathbb R}^n)\mid \operatorname{supp}u\subset \GGcomega \},
\end{equation*}
closed subspace of $H^s({\mathbb R}^n)$; it identifies with the antidual
space of $H^{-s}(\Omega )$ (the space
of antilinear, i.e., conjugate linear, functionals), with a duality consistent with the $L_2$
scalar product. For $s$ equal to a nonnegative integer $k$,
$H^k_0(\GGcomega)$ is usually written $H^k_0(\Omega)$. Spaces over the boundary,
$H^s(\Sigma )$, are defined by local coordinates from $H^s({\mathbb
R}^{n-1})$, $s\in{\mathbb R}$. (There are many equally justified
equivalent choices of
norms there; one can choose a particular norm when convenient.) 
When $s>0$, there are dense continuous embeddings
\begin{equation*}
H^{s}(\Sigma )\subset L_2(\Sigma )\subset H^{-s}(\Sigma ),
\end{equation*}
and there is an identification of
$H^{-s}(\Sigma )$ with the antidual space of $H^s(\Sigma )$, such that the
duality $(\varphi ,\psi )_{-s,s}$ coincides with the $L_2(\Sigma
)$-scalar product when the elements lie there. 
Detailed explanations are found in many
books, e.g.\ 
Lions and Magenes [LM68], H\"o{}rmander [H63], [G09]. (There is a difference of notation: For
$k\in\mathbb N$, the spaces
denoted $H^{k+\frac12}_0(\Omega )$ in \cite{GGLM68} 
are not the same as our $H^{k+\frac12}_0(\GGcomega )$ that are
consistent with \cite{GGH63}; they have the best duality and 
interpolation properties.) 

Consider the case where $A$ is defined on a smooth open subset $\Omega
$ of ${\mathbb R}^n$ and has coefficients in $C^\infty (\GGcomega)$, and
assume that $A$ is elliptic on $\GGcomega$. The
results in this case are a model for results under weaker smoothness hypotheses.
One defines the maximal realization $\GGAma$ as the operator acting like
$A$ in the distribution sense with domain
\begin{equation*}
 D(\GGAma)=\{u\in L_2(\Omega )
\mid Au\in L_2(\Omega )\};
\end{equation*}
it is a closed, unbounded operator in $L_2(\Omega )$.
The minimal realization $\GGAmi$ is
defined as the closure of $A$ acting on $C_0^\infty (\Omega )$ (the
compactly supported $C^\infty $-functions on $\Omega $).
When $\Omega $ is bounded, or is unbounded and there are suitable
bounds on the coefficients of $A$, 
\begin{equation*}
D(\GGAmi)= H^{m}_0(\Omega ).
\end{equation*}
The formal adjoint $A'$ of $A$ is the differential operator acting
as follows:
\begin{equation*}
A'u=\sum_{|\alpha |\le m}D^\alpha (\overline a_\alpha (x) u).
\end{equation*}
By definition, $\GGAmi'$ and $\GGAma$ are adjoints of one another (as
unbounded operators in $L_2(\Omega )$). 

The linear operators $\GGwA$ satisfying 
\begin{equation*}
\GGAmi\subset \GGwA\subset \GGAma 
\end{equation*}
are called {\it realizations} of $A$. 

Generally, $\GGAma$ is far from being
injective, whereas $\GGAmi$ is far from being surjective, but
it may be possible to find realizations $\GGwA$ that are bijective from $D(\GGwA)$ to $L_2(\Omega )$.

We see that the theory of distributions (which in this context was 
preceded by the definition of 
differential operators
acting in the {\it weak sense}) allows defining operators
representing the action of $A$ in a generalized sense. Here
invertibility can sometimes be achieved by methods of functional
analysis.
A fundamental example is the Dirichlet problem (where $\gamma
_0u=u|_\Sigma  $)
\begin{equation}
Au=f \text{ in }\Omega, \quad \gamma_0 u=\varphi \text{ on }\Sigma ,\label{GG2.10}
\end{equation}
for a strongly elliptic second-order operator having $\operatorname{Re}(Av,v)\ge
c\|v\|^2_{L_2(\Omega )}$ with $c>0$ for $v\in C_0^\infty(\Omega
)$.
By use of the so-called Lax-Milgram lemma one
could define a realization $A_\gamma $ of $A$ with $D(A_\gamma )\subset
H^1_0(\Omega)$, such that
$A_\gamma \colon D(A_\gamma )\to L_2(\Omega )$ bijectively. (Details
are found in many books, e.g.\ \cite{GGG09}, Ch.\ 12.)

But then the question was: How close is $A_\gamma ^{-1}$ to solving
the problem in a more classical sense? Second-order derivatives have a
meaning on $H^2(\Omega)$, by closure of the definition on
$C^2(\GGcomega)$, so one can ask: 
\begin{itemize}
\item If $f\in
L_2(\Omega )$, is $u\in H^2(\Omega )$?
\item More generally, if $f\in
H^k(\Omega )$ for some $k\in \mathbb N_0$, is $u\in H^{k+2}(\Omega )$?
\end{itemize}

The answer was first found for the behavior of $u$ in the interior of
$\Omega $: Indeed, when $f\in
H^k(\Omega )$,  $u$ is in $H^{k+2}$ over subsets of $\Omega $ with 
positive distance
from the boundary. This is the so-called {\it interior regularity}.

There remained the question of {\it regularity at the boundary}. It
was answered positively in papers by Nirenberg \cite{GGN55} and Ladyzhenskaya
(see the account in \cite{GGL85}).
This was followed up by research on higher-order operators $A$ and more
general boundary conditions $Tu=\varphi$ (possibly vector valued), where results on interior regularity and
regularity at the boundary were established under suitable
conditions. Besides ellipticity of the operator $A$ one needs a
condition on how the boundary condition fits together with $A$. Some
authors called it the ``covering condition'' or the ``complementing 
condition'', but the name ``the
Shapiro-Lopatinski\u\i{} condition'' (after \cite{GGSh53}, \cite{GGL53}) has been more
generally used. It is also customary to call the system
$\{A,T\}$ {\it elliptic} when it holds (this was suggested in
H\"o{}rmander \cite{GGH63}; we return to a motivation in Section \ref{GG4psdbo}).
A fundamental reference in this connection is the paper of Agmon,
Douglis and Nirenberg \cite{GGADN59} that collects and expands the
knowledge on elliptic boundary value problems. An important point of
view was to obtain so-called ``\`a priori estimates'' (estimates of a
Sobolev norm on $u$ by norms on $Au$ and $Tu$ plus a lower order norm on $u$), shown
for smooth functions at first, and extended to the considered solution. 

Important monographs exposing the theories and the various authors' own
contributions were written by  Agmon \cite{GGA65}, Lions and Magenes
\cite{GGLM68}; the latter moreover contains valuable information on the
surrounding literature. The early theory is exposed in  Courant and
Hilbert I--II \cite{GGCH53}, \cite{GGCH62}.

In the works at that time, although the striving to show existence 
of a solution
operator was always in the picture, the emphasis was more on showing
qualitative properties of the unknown function $u$ in terms of
properties of the
given data $f$ and $\varphi$, regardless of whether $u$ could be
described by an operator acting on $\{f,\varphi\}$ or not.

Direct machinery to construct approximate 
solution operators in general came into
the picture with the advent of pseudodifferential methods.

\section {Pseudodifferential operators}  \label{GG3psdo}

One of the few cases where an elliptic differential operator has an
explicit solution operator is the case of $I-\Delta $ on $\mathbb R^n$,
whose action
can be described by use of the Fourier transform $\mathcal F$ as
$(1-\Delta )u=\mathcal F^{-1}\bigl((1+|\xi |^2)\mathcal F u\bigr)$, and whose
solution
operator is\begin{equation*}
\operatorname{Op}\Bigl(\frac1{1+|\xi |^2}\Bigr)u=
\mathcal F^{-1}\Bigl(\frac1{1+|\xi |^2}\mathcal F u\Bigr).
\end{equation*}
A variable-coefficient elliptic differential operator on $\mathbb R^n$ can also be described
by the help of the Fourier transform,
\begin{equation*}A u=\sum_{|\alpha |\le
m}a_\alpha (x)D^\alpha u =\sum_{|\alpha |\le
m}a_\alpha (x)\mathcal F^{-1}(\xi ^\alpha \mathcal F u)= \mathcal F ^{-1}
a(x,\xi)\mathcal F u
=\operatorname{Op}(a(x,\xi ))u,
\end{equation*}
where $a(x,\xi )=\sum_{|\alpha |\le m}a_\alpha (x)\xi ^\alpha $ is
the symbol. But even when the symbol satisfies $a(x,\xi )\ne 0$ for all
$x,\xi $, the operator \begin{equation*}
\operatorname{Op}(a(x,\xi )^{-1})=\mathcal F^{-1} a(x,\xi )^{-1}\mathcal F 
\end{equation*} 
is not an exact inverse. Nevertheless, it is useful in the discussion
of solutions, since one can show that \begin{equation*}
\operatorname{Op}(a(x,\xi ))\operatorname{Op}(a(x,\xi )^{-1})=I+\mathcal R,
\end{equation*}
where the remainder $\mathcal R$ is of  order $-1$ (lifts the exponent of
a Sobolev space by 1). 

 A thorough treatment of $\operatorname{Op}(a(x,\xi )^{-1})$ and suitable generalizations that come
closer to being an inverse of $A$ (such ``almost-inverses'' are called
parametrices) was obtained with the invention of {\it
pseudodifferential operators}, $\psi $do's. Some of the initiators
were Mihlin \cite{GGM48}, Calder\'o{}n and Zygmund \cite{GGCZ57}, 
Seeley \cite{GGSe65},
Kohn and Nirenberg \cite{GGKN65}, H\"o{}rmander \cite{GGH65}, \cite{GGH71}.

General $\psi $do's are defined from general
symbols $p(x,\xi )$ as 
\begin{equation*}
\operatorname{Op}(p(x,\xi ))u= \mathcal F ^{-1}
p(x,\xi)\mathcal Fu=\int_{\mathbb
R^n}e^{ix\cdot\xi }p(x,\xi )\hat u(\xi )\, \GGdd\xi,
\end{equation*}
where we use the notation $\GGdd\xi =(2\pi )^{-n}d\xi $; here 
$p(x,\xi )$ is required to belong to a suitable class of functions.

Not only do the operators make sense
on $\mathbb R^n$ where the Fourier transform acts, they are also given a meaning
on manifolds, by use of coordinate change formulas and cutoff
functions. The theory is not altogether easy; it uses concepts from
distribution theory in a refined way. Moreover, it is not exact but
qualitative in many statements, so it can be something of a challenge to
derive good results from its use. A fine achievement is that it leads
to Fredholm 
operators, when applied to elliptic operators on compact manifolds. Here
one is  just a small step away from having invertible operators; this
can sometimes be achieved by relying on additional knowledge of
the situation.

A so-called ``classical'' $\psi $do is an operator defined from a
symbol that has an asymptotic series expansion in homogeneous terms (a
polyhomogeneous symbol):
\begin{equation*} 
p(x,\xi)\sim\sum_{j\in{\mathbb N_0}}p_{m-j}(x,\xi),\quad 
p_{m-j}(x,t\xi )=t^{m-j}p_{m-j}(x,\xi )\text{ for }|\xi
|\ge 1, \; t\ge 1.
\end{equation*}
It is said to be of order $m$, and $\GGOp(p)$ maps $H^s$ to $H^{s-m}$ for
all $s\in\mathbb R$. $p_m$ is called the principal symbol, and $p$ is
said to be elliptic when $p_m(x,\xi )\ne 0$ for $|\xi |\ge 1$.  Here one has that\begin{equation*}
\GGOp(p)\GGOp(p')=\GGOp(pp')+\mathcal R_1=\GGOp(p_mp'_{m'})+\mathcal R_2,\end{equation*}
where $\mathcal R_1$ and $\mathcal R_2$ are of order $m+m'-1$. In this way, the principal
part {\it dominates the behavior}. When $p$ is elliptic, the principal
part of a parametrix is found as $p_m^{-1}$ (for $|\xi |\ge 1$,
extended smoothly to $\xi \in \GGrn $). Also the notation $p^0$ is used
for the principal symbol.

There exist many different symbol classes with generalizations of the above
properties, designed for particular purposes.

The very attractive feature of classical $\psi $do's is that they form a scale
of operators of {\it all integer orders}, including differential
operators among those of positive order, and including
parametrices and inverses of elliptic differential operators among
those of negative order. Moreover, it is an ``algebra'', in
the sense that the elements by composition (and by addition) lead to
other classical $\psi $do's.

The calculus is explained in the original papers and in several
subsequent books, such as Treves \cite{GGT80}, H\"o{}rmander \cite{GGH85};
a detailed introduction can be also found in Chapters 7--8 of \cite{GGG09}.

The $\psi $do theory gives (after one has done the work to set it up)
an easy proof of {\it interior regularity} of solutions to elliptic problems.

\section{Pseudodifferential
boundary operators} \label{GG4psdbo}

When an elliptic  differential operator $A$ is considered on a subset of $\GGrn$ or on
a manifold with boundary --- let us here for simplicity just consider the case of a
smooth bounded open subset $\Omega $ of $\GGrn$ --- we must impose
boundary conditions to get uniquely solvable problems. Let assume that
we are in a case where the boundary condition $Tu=\varphi $ together
with \eqref{GG2.4} gives a uniquely solvable problem; here $T$ is a {\it trace
operator} mapping functions on $\Omega $ into $M$-tuples of 
functions on $\Sigma =\partial\Omega $. 
We can formulate this in terms of matrices:
\begin{equation}
\begin{pmatrix} A  \\ \quad \\T
\end{pmatrix}\colon C^\infty (\GGcomega)\to \begin{matrix} C^\infty (\GGcomega)\\ \times\\
C^\infty (\Sigma )^M  \end{matrix} \text{ has an inverse }
\begin{pmatrix} R  & K \end{pmatrix} \colon \begin{matrix} C^\infty (\GGcomega)\\ \times
\\C^\infty (\Sigma
)^M 
\end{matrix}
\to  C^\infty (\GGcomega) .\label{GG4.1}
\end{equation}
Here  $K$ is called a {\it Poisson
  operator}; it solves the semi-homogeneous problem
\begin{equation*}
Av=0\text{ in }\Omega ,\quad Tv=\varphi\text{ on }\Sigma  .
\end{equation*}
The operator $R$ solves the other semi-homogeneous problem
\begin{equation*}
Aw=f\text{ in }\Omega ,\quad Tw=0 \text{ on }\Sigma .
\end{equation*}
In a closer analysis of $R$, we can write it as a sum of two terms:
\begin{equation} R=Q_++G,\label{GG4.4}
\end{equation}
where $Q$ is the $\psi $do $A^{-1}$ on ${\mathbb R}^n$,  
$Q_+$ is its {\it truncation}  $r^+Qe^+$ to $\Omega $, and $G$ is a
supplementing operator adapted to the specific boundary condition,
called a {\it singular Green operator} (s.g.o.).
The operator $e^+$ stands for extension by 0 (to functions on ${\mathbb
R}^n$), and the operator $r^+$ stands for restriction to $\Omega $.

The calculus of pseudodifferential boundary operators ($\psi $dbo's) was
initated by Boutet de Monvel \cite{GGB71}, who introduced operator
systems encompassing both the systems $\binom A T$ and their solution
operators $\begin{pmatrix} R & K\end{pmatrix}$. The original presentation is
somewhat brief, and was followed up by extended expositions, in the
detailed book of Rempel and Schulze \cite{GGRS82}, which elaborated the
index theory, and in the paper
\cite{GGG84} which completed some proofs of composition rules (with
new points of view), and showed spectral asymptotic estimates for
singuar Green operators. The book \cite{GGG96}, whose first edition was issued in 1986,
developed a calculus of parameter-dependent $\psi $dbo's, leading to
resolvent and heat operator constructions. The recent book \cite{GGG09}
gives a full introduction to the theory.

In the systematic calculus of
  Boutet de Monvel one consider
systems (called Green operators):
\begin{equation*}\mathcal A=
\begin{pmatrix} P_++G& K\\ \quad \\T&S 
\end{pmatrix} \colon \begin{matrix} C^\infty (\GGcomega)^N\\ \times \\ C^\infty (\Sigma )^M
\end{matrix} \to \begin{matrix} C^\infty (\GGcomega)^{N'}\\ \times \\ C^\infty (\Sigma )^{M'}
\end{matrix},\text{ where}
\end{equation*}
\begin{itemize}
\item $P$ is a ps.d.o.\  on ${\mathbb R}^n$,
satisfying the so-called {\it transmission condition} at $\Sigma $ (always true for
operators stemming from elliptic differential operators);
\item $P_+=r^+Pe^+$ is the truncation  to $\Omega $ 
(the transmission condition
assures that $P_+$ maps $C^\infty (\GGcomega)$ into $C^\infty
(\GGcomega)$);
\item $T$ is a {\it trace operator} from $\Omega $ to $\Sigma $,
$K$ is a {\it Poisson operator} from $\Sigma $ to $\Omega $,
$S$ is a $\psi $do on $\Sigma $;
\item $G$ is a {\it singular Green operator}, e.g.\ of the type $KT$.
\end{itemize}
The composition of two such systems is again a system belonging to the
calculus. 

The operators extend to act on Sobolev spaces. For $T$ and $G$
there is a condition expressing which differential trace operators 
$\gamma _ju=(\partial/\partial
n)^ju|_{\Sigma }$ that enter: $T$ or $G$ is said to be of class $r$
when only $\gamma _j$'s with $j< r$ enter; 
and then they act on $H^s(\Omega )$ for $s>r-\frac12$.
The class 0 case is the case where
they are purely integral operators, well-defined on $L_2(\Omega )$.

All entries can be matrix-formed. They are defined in local coordinates by formulas
involving Fourier transformation and polyhomogeneous symbols. The idea
is as follows: In local coordinates at the boundary, where $\Omega $
and $
\Sigma $ are replaced by ${\GGrnp}=\{x\in\GGrn\mid x_n>0\}$ and ${\mathbb
R}^{n-1}$ (with points $x'=(x_1,\dots,x_{n-1})$), the system has for
each $(x',\xi ')$ a {\it boundary symbol operator} acting in the $x_n$-variable:
\begin{equation}
\aligned
&a(x',\xi ',D_n)=\\
&\begin{pmatrix} p(x',0,\xi ',D_n) +g(x',\xi ',D_n)& k(x',\xi
',D_n)\\ \quad & \quad \\
t(x',\xi ',D_n)&s(x',\xi ')\end{pmatrix}\colon  \begin{matrix} H^m(\GGrp)^N\\ \times \\ {\mathbb C}^M
\end{matrix} \to \begin{matrix} L_2 (\GGrp)^{N'}\\ \times \\ {\mathbb C}^{M'}
\end{matrix}
\endaligned \label{GG4.5}
\end{equation}
Here $m$ is the order of the operator. Each entry in $a$ acts in a
specific way. E.g., when the matrix is $\binom A T$ in \eqref{GG4.1}, the boundary symbol
operator is the {\it model operator} obtained by freezing the coefficients
at $x'$ and replacing derivatives in $D^\alpha _{x'}$ by their
Fourier transforms $(\xi ')^\alpha $ (with respect to $x'\in {\mathbb R}^{n-1}$). The principal boundary symbol
operator $a^0(x',\xi ',D_n)$ is formed of the top order terms. The principal boundary symbol
operator for $\begin{pmatrix} R&K\end{pmatrix}$ is the inverse of the principal
boundary symbol operator for $\binom A T$. (For \eqref{GG4.5}, $g$ and
$t$ must be of class $\le m$.)

From the boundary symbol operator one defines a full operator by
applying the $\psi $do definition in the $x'$-variable, 
\begin{equation*}
\GGOp'(a(x',\xi ',D_n))u=\int e^{ix'\cdot \xi '}a(x',\xi ',D_n)(\mathcal
F_{y'\to \xi '}u (y',x_n))\,\GGdd \xi '.
\end{equation*}

The symbols have
asymptotic series of terms that are homogeneous in $(\xi ',\xi _n)$
(different rules apply to the different ingredients, and we must refer
to the mentioned references for further details).
One then defines $\mathcal A$ to be {\it elliptic}, when
\begin{enumerate}
\item $P$ is elliptic, i.e.\ its principal symbol $p^0(x,\xi )$ 
is invertible at
each $(x,\xi )$ with $|\xi |\ge 1$,
\item the principal boundary symbol operator $a^0(x',\xi ',D_n)$ is invertible at
each $(x',\xi ')$ with $|\xi '|\ge 1$.
\end{enumerate}
For a system $\binom A T$ formed of an elliptic differential operator $A$ and a differential trace operator $T$, 2) is precisely the old
covering/\-complementing/\-Shapiro-Lopatinski\u\i{} condition for  $\{A,T\}$.

In the elliptic case, one can construct a parametrix $\mathcal B^0$ from the inverses
of the symbols in 1)--2), such that $\mathcal A\mathcal B^0-I$ and $\mathcal
B^0\mathcal A-I$ have order $\le -1$, and the construction can be refined
to give errors of arbitrarily low order. With supplementing
information it can be possible to obtain an inverse.

For example, there holds a a solvability theorem for an elliptic 
differential operator
problem as in \eqref{GG4.1}, formulated in this framework as follows:

\begin{theorem}\label{GGTheorem 4.1} Let $\Omega \subset {\mathbb R}^n$ be a smooth,
bounded open set, denote $\partial\Omega =\Sigma $, and let
$A=\sum_{|\alpha |\le 2m}a_\alpha (x)D^\alpha $ with $a_\alpha
\in C^\infty (\GGcomega)$ be elliptic on
$\GGcomega$, i.e., $\sum_{|\alpha |=2m}a_\alpha (x)\xi ^\alpha \ne 0$
for $x\in\GGcomega$, $ \xi \in {\mathbb R}^n\setminus\{0\}$. 
Let $T=(T_j)_{j=1}^m$ be a column vector of
trace operators $T_j=\gamma _0B_j$, where the $B_j$ are differential
operators of order $m_j$ with $C^\infty $-coefficients, $0\le
m_1<\dots <m_m\le 2m-1$. (Then $T$ is of class $r=m_m+1\le 2m$.) Assume
that $\{A,T\}$ is elliptic.

The operator  $\mathcal A=\binom A T$ defines a continuous mapping
\begin{equation}
\mathcal A=\begin{pmatrix} A  \\ \quad \\T
\end{pmatrix}\colon H^{2m+s}(\Omega)\to \begin{matrix} H^s (\Omega)\\ \times\\
\prod _{j=1}^m H^{2m+s-m_j-\frac12}(\Sigma )\end{matrix}  \text{ for }s>r-2m-\tfrac12,
 \label{GG4.6}
\end{equation} 
and there is a system $\mathcal B=\begin{pmatrix} R& G\end{pmatrix}$ (a parametrix)
belonging to the calculus and continuous in the opposite direction, such that\begin{equation*}\aligned
\mathcal A\mathcal B&=\begin{pmatrix} I&0\\0&I\end{pmatrix} +\mathcal R_1,\quad \mathcal B\mathcal
A=I+\mathcal R_2,\\
\mathcal R_1&\colon \begin{matrix} H^s (\Omega)\\ \times\\
\prod _{j=1}^m H^{2m+s-m_j-\frac12}(\Sigma )\end{matrix} \to \begin{matrix} H^{s'} (\Omega)\\ \times\\
\prod _{j=1}^m H^{2m+s'-m_j-\frac12}(\Sigma )\end{matrix}\\
\mathcal R_2 &\colon  H^{2m+s}(\Omega)\to H^{2m+s'}(\Omega)
\endaligned
\end{equation*}
for all $s$ as in {\rm \eqref{GG4.6}}, all $s'\ge s$. Here $K$ is a row vector
of Poisson operators $(K_j)_{j=1}^m$ of orders $-m_j$, and $R=Q_++G$,
where $Q$ is a parametrix of $A$ on a neighborhood of $\GGcomega$, and
$G$ is a singular Green operator.

The operator $\mathcal A$ in {\rm \eqref{GG4.6}} is Fredholm for each $s$, with
the same finite dimensional  kernel and cokernel in $C^\infty $ for all
$s$.  

If $\mathcal A$ is bijective, the inverse belongs to the calculus (and
is of the same form as $\mathcal B$).

\end{theorem}

When $r=2m$, the lower limit for $s$ is  $ -\frac12$; cases where it is
$<-\frac12$ occur for example for the Dirichlet problem, where $r=m$, and
$s$ can go down to $-m-\frac12$. It is
useful to know that the Poisson operator $K$ in fact has the mapping property 
\begin{equation*}
K\colon \prod _{j=1}^m H^{2m+s-m_j-\frac12}(\Sigma )\to H^{2m+s}(\Omega)
\end{equation*}
for {\it all} $s\in{\mathbb R}$. The trace operator $T$ is called {\it
normal}, when $\gamma _0B_j=b_j\gamma _{m_j}+\sum_{k<m_j}B_{jk}\gamma
_k$ with an invertible coefficient $b_j$ for each $j$. (More general
normal boundary value problems are described below in Section \ref{GG9appl2}.)

For example, for a second-order strongly elliptic
operator with a Dirichlet condition, the operator in the theorem maps as follows:  
\begin{equation*}
\mathcal A=\begin{pmatrix} A  \\ \quad \\\gamma _0
\end{pmatrix}\colon H^{2+s}(\Omega)\to \begin{matrix} H^s (\Omega)\\ \times\\
H^{\frac32+s}(\Sigma )\end{matrix}  \text{ for }s>-\tfrac32,
\end{equation*} 
with parametrices and solution operators continuous in the opposite direction.

Elliptic operators $A$ of odd order
occur mainly as square matrix-formed operators, and there is a similar theorem
for such cases, where also the $B_j$ can be matrix-formed.
Operators of Dirac-type are a  prominent first-order example. 
The matrix
case is also interesting for even-order operators. The
results can moreover be worked out for operators defined on manifolds,
acting in vector bundles.  (See e.g.\ \cite{GGG74}, on the even-order case, for notation and the
appropriate definition of normal boundary conditions.)

\section{Extension theories}\label{GG5ext}

We shall now recall some elements of the functional analysis theory of
extensions of given operators.
This has       
a long history, with prominent contributions from J.\ von Neumann
in 1929 \cite{GGN29}, K.\ Friedrichs in 1934 \cite{GGF34}, M.\ G.\ Kre\u\i{}n
in 1947 \cite{GGK47},
M.\ I.\ Vishik in 1952 
\cite{GGV52}, M.\ S.\ Birman in 1956 \cite{GGB56}, and others. The present author made a
number of contributions in 1968--74 \cite{GGG68}--\cite{GGG74}, completing the
preceding theories and working out
applications to elliptic boundary value problems; further developments
are found e.g.\ in \cite{GGG83}, \cite{GGG84}, and in recent works.

At the same time there was another, separate development of 
abstract extension
theories, where the operator concept gradually began to be 
replaced by the
concept of relations. This development has been
aimed primarily towards applications to ODE, however including 
operator-valued such equations and
Schr\"o{}dinger operators on ${\mathbb R}^n$; keywords in
this connection are:
boundary triples theory, Weyl-Titchmarsh $m$-functions and Kre\u\i{}n
resolvent formulas. Cf.\ e.g.\ Ko\v{c}ube\u{\i} \cite{GGK75}, 
Vainerman \cite{GGV80},
Lyantze and  Storozh \cite{GGLS83}, 
Gorbachuk
and Gorbachuk \cite{GGGG91}, Derkach and Malamud \cite{GGDM91}, Arlinskii \cite{GGA99}, Malamud and Mogilevskii \cite{GGMM02}, 
Br\"uning,
Geyler and Pankrashkin \cite{GGBGP06}, and their references. In later years
there have also been applications to elliptic boundary value problems,
cf.\ e.g.\ Amrein and Pearson \cite{GGAP04}, N.D.~Kopachevski{\u{\i}}
and S.G.~Kre\u{\i}n \cite{GGKK04},
Behrndt and Langer \cite{GGBL07}, Ryzhov \cite{GGR07}, Brown, Marletta, Naboko and
Wood \cite{GGBMNW08}, Gesztesy and Mitrea \cite{GGGM08}, and their references. 

The connection between the two lines of extension theories has been
clarified in a recent work of Brown, Grubb and Wood \cite{GGBGW09}.

At this point we should also mention the recent efforts for problems on
nonsmooth domains: Posilicano and Raimondi \cite{GGPR09}, Grubb
\cite{GGG08}, Gesztesy and
Mitrea \cite{GGGM08,GGGM11}, Abels, Grubb and Wood \cite{GGAGW11}; here
\cite{GGG08,GGGM11,GGAGW11} use  \cite{GGG68}.

In the following, we shall use the notation from
\cite{GGG68}--\cite{GGG74} and \cite{GGBGW09}.

Let there be given a  pair $A_{\min}$, $A'_{\min}$ of closed, densely defined
operators in a Hilbert
space $H$, such that the following holds:
\begin{equation*}
A_{\min}\subset (A'_{\min})^*=:A_{\max},\quad A'_{\min}\subset
(A_{\min})^*=:A'_{\max}.
\end{equation*}
 Let $\mathcal M=\{\GGwA\mid \GGAmi\subset
\GGwA\subset \GGAma\}$. Write $\GGwA u$ as $Au$, when $\GGwA\in\mathcal M$.

We assume that there is given an operator $A_\gamma \in\mathcal M$, the reference
operator, with $0\in
\varrho (A_\gamma )$ (the resolvent set); then
\begin{equation*}A_{\min}\subset A_\gamma \subset A_{\max},\quad
A'_{\min}\subset A_\gamma ^* \subset A'_{\max}.\end{equation*}

The case where $\GGAmi=\GGAmi'$ and $A_\gamma $ is selfadjoint, is called
the {\it symmetric} case.

Let
$Z=\operatorname{ker}\GGAma$, $Z'=\operatorname{ker}\GGAma'$, and
define the basic non-orthogonal decompositions
\begin{equation*}
D(\GGAma)=D(A_\gamma )\dot+ Z,\quad D(\GGAma')=D(A_\gamma ^*)\dot+
Z',
\end{equation*}
 denoted $u=u_\gamma +u_\zeta
=\GGpr_\gamma u+\GGpr_\zeta u$, where $\GGpr_\gamma =A_\gamma ^{-1}\GGAma$,
with a similar notation with primes.

  By $\GGpr_Xu=u_X$ we denote the {\it orthogonal projection} from $H$ to
$X$. The injection  $X\hookrightarrow H$ is denoted $\GGinj_X$ (it is
the adjoint of $\GGpr_X\colon H\to X$).

There holds an ``abstract Green's formula'' for $u\in D(\GGAma)$, $v\in D(\GGAma')$:
\begin{equation}
(A u,v) - (u, A' v)=((A u)_{Z'}, v_{\zeta '})-(u_\zeta ,
(A 'v)_Z).\label{GG5.3}
\end{equation}

It can be used to show that when $\GGwA\in \mathcal M$, and we define
\begin{equation*}
 V=\overline{\GGpr_{\zeta } D(\widetilde A)},\quad W=\overline{\GGpr_{\zeta '} D(\widetilde A^*)},
\end{equation*}
then
\begin{equation*}
\{\{u_\zeta , (A u)_W\}\mid u\in D(\GGwA)\}\text{ is a graph,}
\end{equation*}
defining an operator $T$ from
$D(T)\subset V$ to $W$.

\begin{theorem}\label{GGTheorem 5.1}
(\cite{GGG68}) 
There is a {\rm 1--1} correspondence between the
   closed operators $\GGwA\in {\mathcal M}$ and the closed densely defined
operators $T\colon V\to W$,
   where $V\subset Z$, $W\subset Z'$ (arbitrary closed subspaces), such that
   $\GGwA$ corresponds to $T\colon V\to W$ if and only if
\begin{equation}
D(\GGwA)=\{u\in D(\GGAma)\mid \GGpr_\zeta u\in D(T),\;  (Au)_W= T\GGpr_\zeta u\}.\label{GG5.6}
\end{equation}

In this correspondence, $ V=\overline{\GGpr_{\zeta } D(\widetilde A)}$,
$W=\overline{\GGpr_{\zeta '} D(\widetilde A^*)}$, and
\begin{itemize}
\item $\GGwA^*$ corresponds analogously to $T^*\colon W\to V$.

\item $\operatorname{ker}\GGwA=\operatorname{ker}T$; 
\quad $\operatorname{ran}\GGwA=\operatorname{ran}T+(H\ominus W)$.
 
\item  $\GGwA$ is bijective if and only if $T$ is so, and then 
\begin{equation*}\GGwA^{-1}=A_\gamma
^{-1}+\GGinj_{V }T^{-1}\GGpr_W.
\end{equation*}
\end{itemize}

One also has
\begin{equation*}
D(\GGwA)=\{u=v+A_\gamma ^{-1}(Tz+f)+z\mid v\in D(\GGAmi), \, z\in D(T),\,
f\in Z\ominus W \},
\end{equation*} 
where $v$, $z$ and $f$ are uniquely determined from $u$.
\end{theorem}
 
The result builds on the works of  Kre\u\i{}n \cite{GGK47} and Birman \cite{GGB56} (for
selfadjoint operators) and Vishik \cite{GGV52}, and completes the latter: 
In Vishik's paper, the
$\GGwA$ were set in relation to operators over the nullspaces going in the
opposite direction of our $T$'s, and the results were focused on those
$\GGwA$'s that have closed range (the so-called normally solvable realizations).  Our analysis covered all closed $\GGwA$.

The condition in \eqref{GG5.6}
\begin{equation}(Au)_W= T\GGpr_\zeta u\label{GG5.9}\end{equation} 
can be viewed as an ``abstract boundary condition''. 

When $\lambda \in \varrho (A_\gamma )$, one can do
the same construction for the operators shifted by subtraction of $\lambda $. We denote
 \begin{equation*}
Z_\lambda =\operatorname{ker}(\GGAma-\lambda ), \quad
   Z'_{\bar\lambda }=\operatorname{ker}(\GGAma'-\bar\lambda ),
\end{equation*}
and have the decompositions (where $\GGpr^\lambda _\gamma = (A_\gamma
-\lambda )^{-1}(\GGAma -\lambda )$)
\begin{equation*}D(\GGAma)=D(A_\gamma )\dot+ Z_\lambda ,  \quad u=u^\lambda _\gamma +u^\lambda _\zeta
=\GGpr^\lambda _\gamma u+\GGpr^\lambda _\zeta u,
\end{equation*} 
with a similar notation with primes.

\begin{corollary}\label{GGCorollary 5.2}
 Let $\lambda \in \varrho (A_\gamma )$. For the
closed $\GGwA\in\mathcal M$,  there is a 1--1 correspondence
\begin{equation*}
\GGwA -\lambda  \longleftrightarrow \begin{cases} T^\lambda \colon V_\lambda \to W_{\bar\lambda },\text{ closed, densely
defined}\\
\text{with }V_\lambda \subset Z_\lambda ,\; W_{\bar\lambda }\subset Z'_{\bar\lambda },\text{ closed subspaces.} \end{cases}
\end{equation*}

  Here $D(T^\lambda )=\GGpr^\lambda _\zeta  D(\widetilde A)$, $V_\lambda =\overline{D(T^\lambda )}$, $
W_{\bar\lambda }=\overline{\GGpr^{\bar\lambda }_{\zeta '} D(\widetilde
A^*)}$, and $D(\GGwA)$ consists of the functions $u\in D(\GGAma)$ such
that $u^\lambda _\zeta\in D(T^\lambda )$ and
\begin{equation*}
T^\lambda u^\lambda _\zeta=((A-\lambda )u)_{W_{\bar\lambda }}.
\end{equation*}
   
Moreover,
\begin{itemize}
\item
$\operatorname{ker}(\GGwA-\lambda )=\operatorname{ker}T^\lambda $;
\quad $\operatorname{ran}(\GGwA-\lambda )=\operatorname{ran}T^\lambda +(H\ominus W_{\bar\lambda })$.
 \item
$\GGwA-\lambda $ is bijective if and only if $T^\lambda $ is so, and when $\lambda \in \varrho (\GGwA)\cap \varrho (A_\gamma )$, 
\begin{equation*}(\GGwA-\lambda
)^{-1}=(A_\gamma -\lambda )
^{-1}+\GGinj_{V_\lambda  }(T^\lambda )^{-1}\GGpr_{W_{\bar\lambda
  }}.
\end{equation*}
\end{itemize}

\end{corollary}
 
 This gives a Kre\u\i{}n-type resolvent formula for any closed
$\GGwA\in\mathcal M$ with $\varrho (\GGwA)\cap \varrho (A_\gamma )\ne \emptyset$.

The relation between $T$ and $T^\lambda $ was determined in \cite{GGG74}
in the symmetric case, for real $\lambda $, and the proof given there
extends immediately to the general situation (as shown in \cite{GGBGW09}):
 
For $\lambda \in \varrho (A_\gamma )$, define
\begin{equation*}
E^\lambda =I+\lambda (A_\gamma -\lambda )^{-1},\text{ it has the inverse
}
F^\lambda =I-\lambda A_\gamma ^{-1},
\end{equation*}
and similarly
$E^{\prime\bar\lambda }=I+\bar\lambda (A^*_\gamma -\bar\lambda
)^{-1}$ has the inverse $
F^{\prime\bar\lambda }=I-\bar\lambda (A^*_\gamma )^{-1}$ on $H$.
Then $E^\lambda F^\lambda =F^\lambda E^\lambda =I$,
$E^{\prime\bar\lambda }F^{\prime\bar\lambda }=F^{\prime\bar\lambda
}E^{\prime\bar\lambda }=I$ on $H$.
  
Moreover, $E^\lambda $ and
$E^{\prime\bar\lambda }$ restrict to 
homeomorphisms 
\begin{equation*}
E^\lambda _V\colon  V\GGsimto V_\lambda ,\quad E^{\prime\bar\lambda
}_W\colon W\GGsimto W_{\bar\lambda },
\end{equation*} 
  with inverses $F^\lambda _V$ resp.\ $F^{\prime\bar\lambda }_W$.
In particular, $D(T^\lambda )=E^\lambda _V D(T)$.

The operator families derived from $E^\lambda  $ are related to what was  called a
gamma-field in other works from the 1970's and onwards, as a simple special case.

\begin{theorem}\label{GGTheorem 5.3}
 Let $G^\lambda _{V,W}=-\GGpr_W\lambda
E^\lambda \GGinj_{V }$; then
\begin{equation*}(E^{\prime\bar\lambda }_W)^* T^\lambda E^\lambda _V=T+G^\lambda
_{V,W}.
\end{equation*} 
 
  In other words, $T$ and $T^\lambda $ are related by the
commutative diagram
\begin{equation*}
\CD
V_\lambda     @<\sim<{E^\lambda _{V}} <  V    \\
@V  T^\lambda VV           @VV  T+G^\lambda _{V,W} V\\
W_{\bar\lambda }   @>\sim >(E^{\prime\bar\lambda }_W)^* >   W
\endCD \hskip2cm D(T^\lambda )=E^\lambda _V D(T).
\end{equation*}

\end{theorem}

%
In the joint work with Brown and Wood \cite{GGBGW09} we moreover showed
how the study relates to studies
of boundary triples and $M$-functions by 
other researchers (as referred to in the start af this section;  more references are given in
\cite{GGBGW09}):

Let $V=Z$,
$W=Z'$, then, with ${\mathcal H}=Z'$, ${\mathcal K}=Z$, the mappings
\begin{equation*}\begin{pmatrix} \Gamma _1u\\ \Gamma _0u\end{pmatrix}=\begin{pmatrix} (A u)_{Z'}\\
u_\zeta \end{pmatrix} \colon  D(\GGAma)\to {\mathcal H}\times {\mathcal K},\end{equation*}
\begin{equation*}\begin{pmatrix} \Gamma '_1v\\ \Gamma '_0v\end{pmatrix} =\begin{pmatrix} (A' v)_{Z}\\
v_{\zeta '}\end{pmatrix} \colon  D(\GGAma')\to {\mathcal K}\times 
{\mathcal H},
\end{equation*}
form a boundary triple: Both mappings $\binom {\Gamma _1}{\Gamma _0}$
and $\binom {\Gamma '_1}{\Gamma '_0}$ are surjective,  their
kernels are  $D(\GGAmi)$ resp.\ $D(\GGAmi')$, and they satisfy the Green's formula
\begin{equation*}
(\GGAma u,v) - (u, \GGAma' v)=(\Gamma _1u, \Gamma '_0v)_{\mathcal H}-(\Gamma _0u ,
\Gamma '_1v)_{\mathcal K},
\end{equation*}
which is a rewriting of \eqref{GG5.3}.

Here one can consider a boundary condition 
\begin{equation}
\Gamma _1u=T\Gamma _0u,\label{GG5.20}
\end{equation}
where we allow $T$ to be unbounded (closed densely defined) from $\mathcal
K$ to $\mathcal H$; it defines a restriction $\GGwA$ of $\GGAma$ by
$D(\GGwA)=\{u\in D(\GGAma)\mid \Gamma _0u\in D(T), \Gamma _1u=T\Gamma
_0u\}$. Then it is customary to define an $M$-function as follows:

\begin{definition}\label{GGDefinition 5.4}  For $\lambda \in \varrho (\GGwA)$, $M(\lambda
)\colon \operatorname{ran}(\Gamma _1-T\Gamma _0)\to \mathcal K$ is the operator satisfying
\begin{equation*}
M(\lambda )(\Gamma _1x-T\Gamma _0x)=\Gamma _0x,\text{ for all }x\in
\operatorname{ker}(\GGAma-\lambda )=Z_\lambda .
\end{equation*} 
\end{definition}

The analysis in \cite{GGBGW09} showed that $M(\lambda )$ is a
holomorphic family of operators in $\mathcal L(\mathcal H,\mathcal
K)$. On the other hand, when $\GGwA$ and its boundary condition \eqref{GG5.20} are considered from the point of view of 
 extensions \cite{GGG68}--\cite{GGG74} recalled further above,
$\GGwA$ is the
operator corresponding to $T\colon Z\to Z'$ by Theorem \ref{GGTheorem
  5.1}. Then we find moreover, in terms of the $\lambda $-dependent families
introduced in that context:
\begin{equation*}
M(\lambda )=-(T+G^\lambda _{Z,Z'})^{-1}=-F^\lambda _Z(T^\lambda )^{-1}(F^{\prime\bar\lambda }_{Z'})^*,\text{ when }\lambda \in \varrho (\GGwA)\cap
\varrho (A_\gamma ).
\end{equation*}
This gives the Kre\u\i{}n resolvent formula in the form
\begin{equation*}
(\GGwA-\lambda
)^{-1}=(A_\gamma -\lambda )
^{-1}-\GGinj_{Z_\lambda  }E^\lambda _ZM(\lambda )(E^{\prime\bar\lambda
}_{Z'})^*\GGpr_{Z'_{\bar\lambda }}.
\end{equation*}

For the case of general $V$ and $W$, we could likewise construct an
$M$-function from $W$ to $V$ for $\lambda \in \varrho (\GGwA)$, and
establish a Kre\u\i{}n resolvent formula around it. The following
result is shown in  \cite{GGBGW09}:

\begin{theorem}\label{GGTheorem 5.5}
  Let $\GGwA$ correspond to $T\colon V\to W$. For 
  $\lambda \in \varrho (\GGwA)$ there is a well-defined holomorphic family
$M(\lambda )\in {\mathcal L}(W,V)$:
\begin{equation*}
M(\lambda )=\GGpr_\zeta \big(I-(\GGwA-\lambda )^{-1}(\GGAma-\lambda )\big)A_\gamma
^{-1}\GGinj_{W }.
\end{equation*}
 
When $\lambda \in \varrho (\GGwA)\cap
\varrho (A_\gamma )$, then
\begin{equation*}
M(\lambda )=-(T+G^\lambda _{V,W})^{-1}=-F^\lambda
  _V(T^\lambda )^{-1}(F^{\prime\bar\lambda }_W)^*,
\end{equation*}
and
\begin{equation*}
(\GGwA-\lambda
)^{-1}=(A_\gamma -\lambda )
^{-1}-\GGinj_{V_\lambda  }E^\lambda _VM(\lambda )(E^{\prime\bar\lambda
}_W)^*\GGpr_{W_{\bar\lambda }}.
\end{equation*}

\end{theorem}

To have both $T^\lambda $ (for $\lambda \in \varrho (A_\gamma )$) 
and $M(\lambda )$ (for $\lambda \in\varrho (\GGwA)$) available is an advantage,
since
$\operatorname{ker}(\GGwA-\lambda )=\operatorname{ker}T^\lambda $ and
$\operatorname{ran}(\GGwA-\lambda )=\operatorname{ran}T^\lambda
+(H\ominus W_{\bar\lambda })$ give straightforward eigenvalue
information at the poles of $M(\lambda)$ in $\varrho(A_\gamma)$.

\begin{remark}\label{GGRemark 5.5} The name $M$-function is consistent with the
notation in some
papers that \cite{GGBGW09} refers to, but possibly diverges from others (one
could also use the longer name Weyl-Titchmarsh function). There is a
recent publication of Malamud \cite{GGM10} that exposes related
resolvent formulas
on the basis of \cite{GGMM02}. (Let us remark that \cite{GGM10} seemingly
ascribes a hypothesis of separate surjectiveness of $\Gamma _0$ and
$\Gamma _1$ to the presentation in \cite{GGBGW09}; this is not so.)
\end{remark}

\section{Implementation of the abstract set-up for elliptic
operators}\label{GG6impl}

We shall now recall how the abstract theory is applied to  a concrete
choice of elliptic operator $A$. Here $\GGAma$ and $\GGAmi$ are defined 
as in Section \ref{GG2bvp}; they are closed operators in $H=L_2(\Omega )$. In \cite{GGG68} general even-order
operators were considered, and the reference operator (called 
$A_\gamma$ in Section \ref{GG6impl}) was taken to represent a general normal boundary
condition. To simplify our explanation, we shall here just consider a
second-order strongly elliptic operator $A$ and let $A_\gamma $ stand
for the Dirichlet realization, mentioned after \eqref{GG2.10}. We have by
elliptic regularity that $D(A_\gamma)=H^2(\Omega )\cap H^1_0(\Omega
)$, and we can assume that (a constant has been added to $A$
such that) the lower bound $m(A_\gamma)$ is positive. The
lower bound $m(P)$ of a operator $P$ is defined by
\begin{equation}
m(P)=\inf\{\operatorname{Re}(Pu,u)\mid u\in D(P),\, \|u\|=1\}\ge
-\infty ;\label{GG6.1}
\end{equation}
when it is finite, $P$ is said to be lower bounded.

The trace operator $\gamma_0$ defines a continuous mapping
$H^s(\Omega)\to H^{s-\frac12}(\Sigma)$ for $s>-\frac12$. We shall also
need
a more advanced fact, namely that, as shown by Lions and Magenes (see
e.g.\ \cite{GGLM68}), $\gamma_0$ extends to a mapping $D(\GGAma) \to
H^{-\frac12}(\Sigma )$, and defines {\it homeomorphisms}
\begin{equation*}
\gamma _Z\colon Z\GGsimto H^{-\frac12}(\Sigma ),\quad \gamma _{Z'}\colon Z'\GGsimto H^{-\frac12}(\Sigma ),\label{GG6.2}
\end{equation*}
where $Z$ and $Z'$ are the nullspaces of $\GGAma$ and $\GGAma'$ (not
contained in $H^s(\Omega)$ for $s>0$).
The inverse of $\gamma _Z$ is consistent with the  Poisson  operator $K_\gamma $
solving the semi-homogeneous Dirichlet problems \eqref{GG2.10} with
$f=0$, in the sense that\begin{equation*}
K_\gamma =\GGinj_{Z }\gamma _Z^{-1}.
\end{equation*}
Similarly, the inverse $\gamma _{Z'}^{-1}$ is consistent with the
Poisson solution operator $K'_\gamma $ solving the Dirichlet problem
for $A'$ with $f=0$, $K'_\gamma =\GGinj_{Z' }\gamma _{Z'}^{-1}$. Moreover, 
with $\lambda $-dependence,
\begin{equation*}
K^\lambda _\gamma =\GGinj_{Z_\lambda  }\gamma _{Z_\lambda }^{-1},\quad
K^{\prime\bar\lambda }_\gamma =\GGinj_{Z_{\bar\lambda }' }\gamma
_{Z_{\bar\lambda }'}^{-1},
\end{equation*}
solve the semi-homogeneous Dirichlet problems for $A-\lambda $,
$A'-\bar\lambda $, when $\lambda\in \varrho(A_\gamma)$.

For a closed subspace
$V\subset Z$, let $X=\gamma _0 V\subset H^{-\frac12}(\Sigma )$. Here we denote the
restriction of $\gamma _0$ by $\gamma _V$;
\begin{equation}
\gamma _V\colon V\GGsimto X,\label{GG6.4}
\end{equation}
with a similar notation for $Y=\gamma _0W$ and $\lambda $-dependent cases.
The map $\gamma _V\colon V\GGsimto X$ has the adjoint $\gamma _V^*\colon X^*\GGsimto
V$.
Here $X^*$ denotes the antidual space of $X$, with a duality
coinciding with the scalar product in $L_2(\Sigma )$ when
applied to elements that come from 
$ L_2(\Sigma )$. The duality is written $(\psi ,\varphi )_{X^*,X}$.

We denote \begin{equation}
K_{\gamma
  ,X}= \GGinj_{V }\gamma _V^{-1}\colon X\to V\subset H;\label{GG6.5}
\end{equation} 
it is a Poisson
operator when $X=H^{-\frac12}(\Sigma )$.

Now a given $T\colon V\to W$ is carried over to a closed, densely 
defined operator $L\colon  X\to
Y^*$ by the  definition
\begin{equation*}
L=(\gamma _W^{-1})^*T\gamma _V^{-1},\quad D(L)=\gamma _V D(T);
\end{equation*}
it is expressed in the diagram
\begin{equation*}
\CD
V     @>\sim>  \gamma _V   >    X\\
@VTVV           @VV  L  V\\
   W  @> \sim > (\gamma _W^{-1})^*>   Y^*\endCD 
 \end{equation*}
There is a similar
definition in the $\lambda $-dependent case.

Before formulating the results in a theorem, we shall interpret the
abstract boundary condition \eqref{GG5.9} defining the
realization $\GGwA$, as a concrete condition in terms of $L$.

$A$ has a Green's formula (for sufficiently smooth $u,v$)
\begin{equation}(Au,v)_{\Omega }-(u,A'v)_{\Omega }=(\nu _1u,\gamma _0v)_\Sigma
-(\gamma _0u, \nu ' _1v)_{\Sigma },\label{GG6.8}\end{equation}
where 
\begin{equation*}\nu _1=s\gamma _1,\quad \nu '_1=\bar s\gamma
_1+\mathcal A'\gamma _0,
\end{equation*}
with a nonvanishing smooth function $s$ and a suitable first-order
differential operator $ \mathcal A'$ on $\Sigma $.

Let $\lambda \in \varrho (A_\gamma )$. In addition to the Poisson operators $K_\gamma ^\lambda $ resp.\
$K^{\prime \bar\lambda }_\gamma $ solving the
Dirichlet problems for $A-\lambda $ resp.\ $A'-\bar\lambda $, 
we shall need the {\it Dirichlet-to-Neumann} operators
\begin{equation*}
P^\lambda _{\gamma _0,\nu _1}=\nu _1K^\lambda _\gamma ,
\quad P^{\prime\bar\lambda }_{\gamma _0,\nu '_1}=\nu '_1
K^{\prime\bar\lambda }_\gamma ,
\end{equation*}
that map the Dirichlet boundary value into the Neumann boundary value
for null-solutions. By the composition rules for $\psi $dbo's, they
are pseudodifferential operators of order 1; moreover, it is known that
$P^\lambda _{\gamma _0,\nu _1}$ is elliptic.

\begin{theorem}\label{GGTheorem 6.1}
 Define the  {\bf reduced Neumann trace operator}
  $\Gamma $ by 
\begin{equation*}\Gamma = 
 \nu _1-P^0_{\gamma _0,\nu  _1}\gamma _0=\nu _1A_\gamma
 ^{-1}\GGAma\colon D(\GGAma)\to  H^{\frac12}(\Sigma ).
\end{equation*}
It is continuous and surjective, and vanishes on $Z$. With the analogous definition for
$A'$ one has the {\bf reduced Green's formula}: 
\begin{equation*}
(Au,v)-(u,A'v)=(\Gamma u,\gamma _0v)_{\frac12,-\frac12}-(\gamma _0u,
\Gamma 'v)_{-\frac12,\frac12}, 
\end{equation*}
valid for all $u\in D(\GGAma)$, $v\in D(\GGAma')$. In particular,
\begin{equation}(Au,w)=(\Gamma u,\gamma _0w)_{\frac12,-\frac12},\text{ when }u\in
D(\GGAma), w\in Z'.\label{GG6.13} 
\end{equation}

For $\lambda \in \varrho (A_\gamma )$ we similarly define
\begin{equation*}
\aligned
\Gamma ^\lambda &= 
 \nu _1-P^\lambda _{\gamma _0,\nu  _1}\gamma _0=\nu _1(A_\gamma -\lambda )
 ^{-1}(\GGAma-\lambda ),\\
\Gamma ^{\prime\bar\lambda} &= 
 \nu _1-P^{\prime\bar\lambda }_{\gamma _0,\nu ' _1}\gamma _0=\nu '_1(A_\gamma ^*-\bar\lambda )
 ^{-1}(\GGAma'-\bar\lambda ),
\endaligned
\end{equation*}
continuous and surjective from $D(\GGAma)$ resp.\ $D(\GGAma')$ to  $H^{\frac12}(\Sigma )$; then
there holds a reduced Green's formula
\begin{equation*}
(Au,v)-(u,A'v)=(\Gamma ^\lambda u,\gamma _0v)_{\frac12,-\frac12}-(\gamma _0u,
\Gamma ^{\prime\bar\lambda }v)_{-\frac12,\frac12}, 
\end{equation*}
 for $u\in
D(\GGAma)$, $ v\in D(\GGAma^*)$. 
\end{theorem}

Here $D(\GGAma)$ is provided with the graph-norm.

Now let $\GGwA$ be a closed operator lying between $\GGAmi$ and $\GGAma$,
so $\GGwA\in\mathcal M$.
The abstract boundary condition \eqref{GG5.9} for $\GGwA$ may be written:
\begin{equation}
(Au,w)=(Tu_\zeta ,w), \text{ all }w\in W.   \label{GG6.16}
\end{equation}
The left-hand side equals $(\Gamma u,\gamma _0
w)_{\frac12,-\frac12}$ by \eqref{GG6.13}. 
The right-hand side equals
\begin{equation*}
(Tu_\zeta ,w)=(T\gamma _V^{-1}\gamma _0u, \gamma _W^{-1}\gamma _0w)=(L\gamma _0u,\gamma _0w)_{Y^*,Y},
\end{equation*}
by definition of $L$ (it is used that when $u_\zeta \in V$,
$u_\zeta=\gamma _V^{-1}\gamma _Vu_\zeta $  $=\gamma _V^{-1}\gamma _0u $).

Hence \eqref{GG6.16} may be written
\begin{equation}
(\Gamma u,\gamma _0
w)_{\frac12,-\frac12}=(L\gamma _0u,\gamma _0w)_{Y^*,Y},\text{ all
}w\in W.\label{GG6.17}
\end{equation}

The injection $\GGinj_Y\colon  Y\to H^{-\frac12}(\Sigma )$ has as adjoint the mapping
$\GGinj _{Y}^*\colon  H^\frac12(\Sigma )\to Y^*$ that sends a functional $\psi $
on $H^{-\frac12}(\Sigma )$ over into a functional $\GGinj_Y^*\psi $ on $Y$ by:
\begin{equation*}
(\GGinj_Y^*\psi ,\varphi )_{Y^*,Y}=(\psi ,\varphi
)_{\frac12,-\frac12}\text{ for all }\varphi \in Y. 
\end{equation*}

Then \eqref{GG6.17} may be rewritten as
\begin{equation*}
\GGinj_Y^*\Gamma u = L\gamma _0u,
\end{equation*}
or, when we use that $\Gamma =\nu _1-P^0_{\gamma _0,\nu _1}\gamma _0$,
\begin{equation}
\GGinj_Y^*\nu _1 u= (L+\GGinj_Y^*P^0_{\gamma _0,\nu _1})\gamma _0u.\label{GG6.20}
\end{equation}
This is the boundary condition derived from \eqref{GG5.9}, when $\GGwA$
corresponds to $T\colon V\to W$ by Theorem \ref{GGTheorem 5.1}, 
carried over to $L\colon X\to Y^*$  by \eqref{GG6.4}. 

Then Theorem \ref{GGTheorem 5.1} 
implies:

\begin{theorem}\label{GGTheorem 6.2}
There is a {\rm 1--1} correspondence between the
   closed operators $\GGwA\in {\mathcal M}$ and the closed densely defined
operators $L\colon X\to Y^*$,
   where $X$ and $Y$ are closed subspaces of $H^{-\frac12}(\Sigma )$, such that
   $\GGwA$ corresponds to $L\colon X\to Y^*$ if and only if
$
D(\GGwA)$ consists of the functions in $ D(\GGAma)$ for which 
\begin{equation}
\gamma _0 u\in D(L),\quad \GGinj_Y^*\nu _1 u= (L+\GGinj_Y^*P^0_{\gamma
  _0,\nu _1})\gamma _0u. \label{GG6.21}
\end{equation}

In this correspondence, $ X=\overline{\gamma _0 D(\widetilde A)}$,
$Y=\overline{\gamma _0 D(\widetilde A^*)}$, and 
\begin{itemize}
\item $\GGwA^*$ corresponds analogously to $L^*\colon Y\to X^*$.

\item $\operatorname{ker}\GGwA=\GGinj_V\gamma _V^{-1}\operatorname{ker}L$; 
\quad $\operatorname{ran}\GGwA=\gamma _W^*\operatorname{ran}L+(H\ominus
W)$, cf.\ {\rm \eqref{GG6.4}, \eqref{GG6.5}}.
 
\item $\GGwA$ is bijective if and only if $L$ is so, and then 
\begin{equation}\GGwA^{-1}=A_\gamma
^{-1}+\GGinj_{V }\gamma _V^{-1}L^{-1}\gamma _W^*\GGpr_W=A_\gamma ^{-1}+K_{\gamma ,X}L^{-1}(K'_{\gamma ,Y})^*.\label{GG6.22}\end{equation}
\end{itemize}
\end{theorem}

Theorems \ref{GGTheorem 6.1} and \ref{GGTheorem 6.2}  are from  \cite{GGG68}, except that we have modified
the notation a little. 

In  \cite{GGBGW09}, the subspace cases are treated with
insertion of an {\it isometry} \linebreak
$\Lambda _\frac12\colon  L_2(\Sigma )\GGsimto H^{-\frac12}(\Sigma )$, that allows replacing
$X$ and $Y$ by closed subspaces $X_1$ and $Y_1$ of $L_2(\Sigma )$,
identified with their dual spaces; then $\GGinj_Y^*$ is replaced by an
orthogonal projection $\GGpr_{Y_1}$.
\smallskip

In the case where $Y=H^{-\frac12}(\Sigma )$, i.e., $W=Z'$, the map $\GGinj_Y^*$
is superfluous, and the second condition in \eqref{GG6.21} takes the form
\begin{equation*}
\nu _1u=(L+P^0_{\gamma _0,\nu _1})\gamma _0u.
\end{equation*}
When also $X=H^{-\frac12}(\Sigma )$, we say that $\GGwA$ represents 
 a {\it Neumann-type condition}
\begin{equation}
\nu _1u=C\gamma _0u; \text{ here }C=L+P^0_{\gamma _0,\nu _1} \text{ on }D(L).\label{GG6.24}
\end{equation}
In this case, $L$ can act like a pseudodifferential
operator, namely when $C$ (in the condition  $\nu _1u=C\gamma
_0u$) is a differential or pseudodifferential operator. 

Let us consider a slightly different set-up where $C$ is a {\it given} first-order differential or
pseudodifferential operator on $\Sigma$, and we {\it define}
$\GGwA$ by
\begin{equation}
D(\GGwA)=\{u\in D(\GGAma)\mid \nu _1u=C\gamma
_0u\},\label{GG6.25}
\end{equation}
where $\gamma _0$ and $\nu_1$ are considered as mappings from
$D(\GGAma)$ to $H^{-\frac12}(\Sigma )$ resp.\ $H^{-\frac32}(\Sigma )$.
We shall discuss the corresponding operator $L\colon X\to Y^*$.
Since \linebreak$\{\gamma_0u,\nu_1u\}\colon H^2(\Omega)\to
H^\frac32(\Sigma)\times H^\frac12(\Sigma)$, $D(L)=\gamma_0D(\GGwA)\supset
H^\frac32(\Sigma)$. Then since $H^\frac32(\Sigma)$ is dense in $H^{-\frac12}(\Sigma)$,
$X=H^{-\frac12}(\Sigma)$. By use of Green's formula \eqref{GG6.8} it is checked
that the adjoint $\GGwA^*$ extends the realization of
$A'$ with domain consisting of the functions $v\in H^2(\Omega)$
satisfying \begin{equation*}
\nu_1'v=C^*\gamma_0 v,
\end{equation*}
so also $Y=H^{-\frac12}(\Sigma)$.
Thus we are in the case of Neumann-type boundary conditions, so by
comparison with \eqref{GG6.2}, it is seen that $L$ acts like
$C-P^0_{\gamma_0,\nu_1}$.

The domain $D(L)$ equals $\{\varphi \in H^{-\frac12}(\Sigma )\mid (C-P^0_{\gamma _0,\nu _1
})\varphi \in H^{\frac12}(\Sigma )\}$; it
may not be easy to determine more exactly. Note that $L$ is used as a  map
from its domain in $H^{-\frac12}(\Sigma)$  to $H^{+\frac12}(\Sigma)$, 
although it acts like a
$\psi$do of order 1.

One case is clear, though: If $C-P^0_{\gamma_0,\nu_1}$ is {\it
elliptic} of order 1, then $L\varphi \in H^{\frac12}(\Sigma)$ implies 
$\varphi \in H^{\frac32}(\Sigma)$; in this case $D(L)=
H^{\frac32}(\Sigma)$. Moreover, $D(\GGwA)\subset
H^2(\Omega )$. A check of the boundary symbol rules shows that
this is precisely the case where the system $\{A,\nu_1-C\gamma_0\}$ is
{\it elliptic}. Here we have:

\begin{theorem}\label{GGTheorem 6.3} Let $C$ be a first-order differential or
pseudodifferential operator on $\Sigma $ and define the realization
$\GGwA$ of $A$ by {\rm \eqref{GG6.25}}. Then if $C-P^0_{\gamma_0,\nu_1}$ is
elliptic, the operator $L\colon X\to Y^*$ corresponding to $\GGwA$ by
Theorem {\rm \ref{GGTheorem 6.2}} acts like $C-P^0_{\gamma_0,\nu_1}$ and has
\begin{equation*}
X=Y=H^{-\frac12}(\Sigma ),\quad D(L)=H^{\frac32}(\Sigma ).
\end{equation*}
Moreover, $D(\GGwA)\subset H^2(\Omega )$. Related statements hold for
the adjoint $\GGwA^*$.

\end{theorem}

A Robin
condition $\nu _1u=b\gamma _0u$, with a smooth function $b(x)$ on
$\Sigma $, is elliptic, since $L$ acts like $b-P^0_{\gamma _0,\nu
  _1}$, where $P^0_{\gamma _0,\nu _1}$ is elliptic of order 1 and $b$ is  of order 0.

In the case of Theorem \ref{GGTheorem 6.3}, 
when $L$ is bijective, the formula \eqref{GG6.22} has the form
\begin{equation}\GGwA^{-1}=A_\gamma ^{-1}+K_{\gamma }L^{-1}(K'_{\gamma
  })^*,\label{GG6.28}
\end{equation}
where all ingredients belong to the $\psi $dbo calculus: $K_{\gamma }$
is a Poisson operator, $L^{-1}$ is a $\psi $do on $\Sigma $, $(K'_{\gamma })^*$ is
a trace operator of class 0, and the composition
$K_{\gamma }L^{-1}(K'_{\gamma })^*$ is a singular Green
operator (of class 0). $A_\gamma ^{-1}$ is the sum $Q_++G$ of a truncated $\psi $do
$Q$ on ${\mathbb R}^n$ and a singular Green operator, as in \eqref{GG4.4}.
 
\medskip

Let us include a few words on higher-order cases:
When $A$ is of order $2m$, there is a Green's formula generalizing
\eqref{GG6.8}, where $\gamma _0$ is replaced by an $m$-vector $\gamma
=\{\gamma _0,\dots,\gamma _{m-1}\}$,  and $\nu _1$ and $\nu '_1$ are replaced by
$m$-vectors of trace operators of orders $m$, \dots, $2m-1$, mapping
into products of Sobolev spaces of different orders over $\Sigma $. One then
gets vector versions of the reduced Neumann trace operators
$\Gamma $ and $\Gamma '$, with matrix-formed versions of the
Dirichlet-to-Neumann pseudodifferential operators, but the basic
mechanisms in the interpretation are the same. There are interesting
cases of 
subspaces $X,Y$ of the products of Sobolev spaces over $\Sigma $, where
ellipticity considerations are relevant. Details are given in
\cite{GGG68}--\cite{GGG74} and \cite{GGBGW09}.

\section{Resolvent formulas} \label{GG7res}

When $\lambda \in
\varrho (A_\gamma )$, there is a similar representation of $\GGwA-\lambda $ in terms of a
boundary condition defined from an operator $L^\lambda $ acting over
the boundary. 
Here it is of particular interest to find the connection
between $L$ and $L^\lambda $, just as we found the connection between
$T$ and $T^\lambda $. 
It turns out that the relation between   $L$ and
$L^\lambda $ is simpler: They both go from $X$ to $Y^*$, whereas $T$
resp.\ $T^\lambda $ map between different spaces due to the shift from
$Z$ to $Z_\lambda $. 
This holds, since 
\begin{equation*}
D(L)=\gamma _0D(\GGwA)=\gamma _0D(\GGwA-\lambda )=D(L^\lambda ),\quad
X=\overline{D(L)}=\overline{D(L^\lambda )}, 
\end{equation*}
with similar statements for $D(L^*)$, $D((L^{\lambda })^*)$ and
$Y$. 
Then we have:

\begin{equation*}
\CD
V     @>\sim>{E^\lambda _{V}} >V _\lambda     @>\sim>  \gamma _{V_\lambda }  >    X\\
@V  T+G^\lambda _{V,W} VV @VT^\lambda VV           @VV  L^\lambda   V\\
 W @>\sim>(F^{\prime\bar\lambda }_{W})^* >  W_{\bar\lambda }
@>\sim>(\gamma _{W_{\bar\lambda }}^*)^{-1} >   Y^*\endCD 
\end{equation*} 
The horizontal maps compose as 
$\gamma _{V_\lambda }E^\lambda _V=\gamma _V$, $(\gamma _{W_{\bar\lambda }}^*)^{-1}(F^{\prime\bar\lambda }_{W})^*=(\gamma _{W}^*)^{-1}$,
so
\begin{equation*}   
L^\lambda =\gamma _V^{-1}(T+G^\lambda _{V,W})(\gamma _{W}^*)^{-1}. 
\end{equation*}

In
terms of $L^\lambda $,
the boundary condition reads (analogously to \eqref{GG6.20}):
\begin{equation}
\GGinj_Y^*\nu _1u=(L^\lambda +\GGinj_Y^*P^\lambda _{\gamma _0,\nu _1})\gamma
_0u,\; \gamma _0u\in D(L).\label{GG7.4}
\end{equation}

Since $D(\GGwA-\lambda )=D(\GGwA)$ is at the same time defined by the boundary condition
$\GGinj_Y^*\nu _1u=(L +\GGinj_Y^*P^0 _{\gamma _0,\nu _1})\gamma _0u$ for $\gamma
_0u\in D(L)$, we have that
$L^\lambda +\GGinj_Y^*P^\lambda _{\gamma _0,\nu _1}=L +\GGinj_Y^*P^0
_{\gamma _0,\nu _1}$ on $D(L)$, so 
\begin{equation*}
 L^\lambda =L+\GGinj_Y^*(P^0 _{\gamma _0,\nu _1}-P^\lambda _{\gamma
   _0,\nu _1})\text{ on }D(L).
\end{equation*}
The last formula is convenient, since $P^0 _{\gamma _0,\nu _1}-P^\lambda _{\gamma
   _0,\nu _1}$ can be shown to be {\it bounded} from
$H^{-\frac12}(\Sigma )$ to $H^\frac12(\Sigma )$; hence
 $L^\lambda $ is a perturbation of $L$ by a bounded operator.

Also the general $M$-function defined in Section \ref{GG5ext} carries over to an
$M$-function on the boundary, a holomorphic family of operators
$M_L(\lambda )\in \mathcal L(Y^*,X)$ defined for $\lambda \in \varrho
(\GGwA)$. 

The results are collected in the following theorem:

\begin{theorem}\label{GGTheorem 7.1}
 Let $\GGwA$ correspond to $T\colon V\to W$ as in Theorem {\rm
   \ref{GGTheorem 5.1}},
carried over to $L\colon X\to Y^*$ as in Theorem {\rm \ref{GGTheorem 6.2}}.
For $\lambda \in\varrho (A_\gamma )$ it is also described by the
boundary condition {\rm \eqref{GG7.4}}, and there holds:

{\rm (i)} For $\lambda \in \varrho (A_\gamma )$,  $P^0 _{\gamma
_0,\nu  _1}-P^\lambda _{\gamma _0,\nu  _1}\in {\mathcal L}(H^{-\frac12}(\Sigma ),
H^{\frac12}(\Sigma ))$ and
\begin{equation*}
L^\lambda =L+\GGinj_Y^*(P^0 _{\gamma _0,\nu  _1}-P^\lambda _{\gamma
  _0,\nu  _1})\text{ on }D(L).
\end{equation*}

{\rm (ii)} For $\lambda \in \varrho (\GGwA)$, there is a related
$M$-function $\in{\mathcal L}(Y^*, X)$,
\begin{equation*}
M_L(\lambda ) =\gamma _0\bigl(I-(\GGwA-\lambda )^{-1}(\GGAma-\lambda
)\bigr)A_\gamma  ^{-1}\GGinj_{W }\gamma _{W}^*.
\end{equation*}

{\rm (iii)} For $\lambda \in \varrho (\GGwA)\cap \varrho (A_\gamma )$,
\begin{equation*}
M_L(\lambda )=-(L^\lambda )^{-1}=-\big(L+\GGinj _{Y}^*(P^0 _{\gamma _0,\nu  _1}-P^\lambda _{\gamma
_0,\nu  _1})\GGinj_{X}\big)^{-1},
\end{equation*}
and we have the
Kre\u\i{}n-type resolvent formulas:
\begin{equation}
(\GGwA-\lambda )^{-1}
-(A_\gamma -\lambda )^{-1}=K^\lambda _{\gamma ,X} 
(L^\lambda )^{-1}(K^{\prime\bar\lambda }_{\gamma ,Y})^*
=-K^\lambda _{\gamma ,X} 
M_L(\lambda )(K^{\prime\bar\lambda }_{\gamma ,Y})^*
.\label{GG7.8}\end{equation}
\end{theorem}
In the case of a Neumann-type boundary condition as in \eqref{GG6.25}, 
$ L^\lambda =C-P^\lambda _{\gamma _0,\nu_1 }$ on $D(L^\lambda )=D(L)$.

\section{Applications of pseudodifferential methods I: \\ Conditions for
lower boundedness} \label{GG8appl1}

The formulas we have shown so far use the terminology of $\psi $dbo's
mainly as a way to indicate what the ingredients in certain
operator compositions are. 
The next question to consider is how properties of $\GGwA$
are reflected in properties of  $L$.
Part of the analysis can be carried out with methods of functional
analysis, but there also exist problems that are solved most efficiently
by involving deeper pseudodifferential principles.

An example of how functional analytic and pseudodifferential
methods are useful together, is the question of lower boundedness inequalities.

We here restrict the attention to the symmetric set-up where $A$ is
formally selfadjoint (so $\GGAma=\GGAmi^*$) and  $A_\gamma$ is
selfadjoint; methods for 
extending the results to nonsymmetric set-ups are found in \cite{GGG74}.
We assume that $A_\gamma $ has positive lower bound $m(A_\gamma )$ (cf.\ \eqref{GG6.1}).

In the symmetric case, the general extensions $\GGwA$ can of course be
nonsymmetric (since $\GGAma$ is so). Let
us speak of the ``selfadjoint case'' when only selfadjoint ${\GGwA}\,$'s
are considered.

In the following, we assume throughout that $\GGwA$ corresponds to
$T\colon V\to W$ as in Theorem \ref{GGTheorem 5.1}, and to $L\colon X\to Y^*$ as in
Theorem \ref{GGTheorem 6.2}. The lower boundedness problem is the problem of how information on lower
bounds on $\GGwA$  is related to similar information on $T$ or $L$.
The following was shown in \cite{GGG70} (together with studies of 
 coerciveness estimates):

\begin{theorem}\label{GGTheorem 8.1} In the symmetric set-up with $A_\gamma$
selfadjoint positive, let $\GGwA$ correspond
to $T\colon V\to W$ as in Theorem {\rm \ref{GGTheorem 5.1}}.

$1^\circ$ If $V\subset W$ and $T$ has a lower bound $m(T)$ satisfying 
$m(T)>-m(A_\gamma)$, then $m(\GGwA)\ge m(T)m(A_\gamma )/(m(T)+m(A_\gamma
))$.

$2^\circ$ Assume that $A_\gamma$ is the Friedrichs extension of
$\GGAmi$.
If $m(\GGwA)>-\infty$, then $V\subset W$ and $m(T)\ge m(\GGwA)$.
\end{theorem}

In the selfadjoint case these rules go back to Birman  \cite{GGB56},
preceded by sesquilinear form results by Kre\u\i{}n \cite{GGK47}.

The properties of $T$ are easily translated to similar properties of $L$ using
the homeomorphism \eqref{GG6.2}; here when $X\subset Y$, we set 
\begin{equation}
m_{-\frac12}(L)=\inf\{\operatorname{Re}(L\varphi ,\varphi )_{Y^*,Y}\mid \varphi \in D(L),\, \|\varphi \|_{-\frac12}=1\},\label{GG8.1}
\end{equation}
for some choice of norm $\|\varphi \|_{-\frac12}$ on $H^{-\frac12}(\Sigma )$
. (One can let $\gamma _Z$ be an isometry, to carry numerical
information over between $T$ and $L$. Sometimes qualitative objects
such as the sign of $m_{-\frac12}(L)$ are sufficiently interesting.) To take $A_\gamma$ as the Friedrichs
extension of $\GGAmi$ means that it is taken as the Dirichlet realization.

In $1^\circ$ we see that the statement ``$m(T)>-\infty \implies
m(\GGwA)>-\infty $'' holds under the additional assumption that
$m(T)>-m(A_\gamma )$; there is a nontrivial question of when that
assumption can be removed.
In \cite{GGG74} it was shown that when $A_\gamma$ is the
Friedrichs extension and $A_\gamma^{-1}$ is a {\it compact} operator, then 
$m(T)>-\infty$ does imply  $m(\GGwA)>-\infty$. This same result was also
announced by Mikhailets and Gorbachuk in \cite{GGGM76} for the
selfadjoint case.

In the application to boundary value problems, we therefore have from
this early result that lower
boundedness of $\GGwA$ and $L$ hold simultaneously when we consider
problems on {\it bounded domains} $\Omega$, for then  $A_\gamma ^{-1}$ is
indeed compact.

For unbounded domains, the question has, to our knowledge,  remained 
unsolved up until recently.  The question is closely connected 
with the comparison of
$T$ with $T^\lambda$ as in Theorem \ref{GGTheorem 5.3}. 
Indeed, as shown in
\cite{GGG74}:

 \begin{proposition}\label{GGProposition 8.2} Let $G^\lambda_{V,W}$ be as defined in
Theorem {\rm \ref{GGTheorem 5.3}}. 
The property
\begin{equation}
m(G^\mu_{Z,Z})\to \infty \text{ for }\mu\to
-\infty,\; \mu\in\mathbb R,
\label{GG8.2}\end{equation}
 is necessary and sufficient for the validity of \begin{equation}
m(T)>-\infty \implies
m(\GGwA)>-\infty\label{GG8.3}
\end{equation}
for general closed $\GGwA\in\mathcal M$.
\end{proposition}

The question was also studied later by Derkach and Malamud \cite{GGDM91}
who worked out an analysis that generalizes Proposition 8.2 and gives further
conditions for the validity of the conclusion from $m(T)$ to $m(\GGwA)$.
However this did not capture elliptic problems on unbounded domains
($n\ge 2$). 

Because of the recent interest in the analysis of extensions, we have
considered the problem again, and found a solution in \cite{GGG11a}
for {\it exterior domains} (complements in ${\mathbb R}^n$ of bounded domains).

\begin{theorem}\label{GGTheorem 8.3} Let $\Omega $ be the complement of a smooth
bounded set $\GGcomega_-$ in $\mathbb R^n$, and let $A$ be symmetric and uniformly
strongly elliptic on $\Omega$ with coefficients in $C^\infty
_b(\GGcomega)$, and with a positive lower bound for $A_\gamma $. In the
application of the extension theory to this situation, {\rm
\eqref{GG8.2}} holds, and hence also {\rm \eqref{GG8.3}}.
\end{theorem}

Here $C^\infty _b(\GGcomega)$ stands for the $C^\infty$-functions that
are bounded with bounded derivatives.

The proof relies on the ``translation'' of abstract operators $\GGwA$ to
concrete operators defined by boundary conditions. Indeed, it turns out
that the lower bound of $G^\mu_{Z,Z}$ behaves like the lower bound
$m_{-\frac12}(Q^\mu )$ (cf.\ \eqref{GG8.1}) of 
$Q^\mu=P^0_{\gamma_0,\nu_1}-P^\mu_{\gamma_0,\nu_1}$. Then the deep part of the proof lies in
setting the operator $Q^\mu$ in relation to the analogous operator
for the interior domain $\Omega_{-}$,
$Q^\mu_{-}$, which does have the desired property, in
view of our knowledge of problems on bounded domains. The point is to
show that $|((Q^\mu -Q^\mu_{-})\varphi ,\varphi )|$
is bounded by $c\|\varphi \|^2_{L_2(\Sigma )}$ uniformly  for $\mu \to -\infty $, so
that addition of $Q^\mu -Q^\mu_{-}$ to $Q^\mu_{-}$ does not violate the
growth of the lower bound. This goes by a delicate application of the
$\psi$dbo calculus. Details are in \cite{GGG11a}.

Let us mention that there is a considerably easier result that holds
regardless of boundedness of $\partial\Omega $ and only requires some
uniformity in the estimates of coefficients of the operators, namely
preservation of coerciveness inequalities
(G\aa{}rding inequalities). We here assume that $\Omega$ is a subset
of $\mathbb R^n$ with smooth boundary, admissible as defined in
\cite{GGG96} (besides bounded domains, this allows exterior domains,
perturbed halfspaces and other cases that can be covered with a finite
system of local coordinates of a relatively simple kind). Moreover, we
assume that 
$A$ is uniformly strongly elliptic on $\overline\Omega
$ with coefficients in $C^\infty _b(\GGcomega)$.
The result 
is that then $\GGwA$ satisfies a G\aa{}rding inequality (with $c>0$, $k\in\mathbb R$) 
\begin{equation*}\operatorname{Re}(\GGwA
  u,u)\ge c\|u\|_1^2-k\|u\|_0^2,\quad u\in D(\GGwA),
\end{equation*}
if and only if $X\subset Y$ and $ L$ satisfies an inequality
\begin{equation}
\operatorname{Re}(L\varphi ,\varphi )\ge c'\|\varphi
\|^2_{ H^\frac12(\Sigma)} - k'\|\varphi \|^2_{
H^{-\frac12}(\Sigma)},\quad \varphi\in D(L).\label{GG8.5}
\end{equation}

In the case of differential or pseudodifferential Neumann-type boundary conditions, 
the inequality \eqref{GG8.5} for $L$ holds
precisely when the pseudodifferential operator it acts like, is
strongly elliptic. Details for bounded sets are in  \cite{GGG70} for
 realizations of $2m$-order operators; the extension
to unbounded sets is shown in \cite{GGG11a} --- the argumentation just
involves standard trace theorems and interpolation inequalities.


\section{Applications of pseudodifferential methods II: 
\\ Spectral asymptotics} \label{GG9appl2} 

In the symmetric set-up, when $\Omega $ is bounded, the eigenvalues of
the selfadjoint operator $A_\gamma  $ form a sequence $\lambda _j$
going to $\infty $ on ${\mathbb R}$. Already in 1912, Weyl \cite{GGW12} showed the
famous estimate  for $A=-\Delta $, for $n=2,3$, $m=2$:
\begin{equation}
\lambda _j(A_\gamma )-c_0j^{m/n}\text{ is }o(j^{m/n})\text{ for }j\to \infty ,\label{GG9.1}
\end{equation}
where $c_0$ is a constant defined from the volume of $\Omega
$; the eigenvalues are repeated according to
multiplicities. Equivalently, the counting function $N(t; A_\gamma )$
(counting the number
of eigenvalues in $[0,t]$), and the eigenvalues $\mu _j$ of the inverse
$A_\gamma ^{-1}$, satisfy
\begin{equation*}
\aligned
N(t;A_\gamma )&-c_At^{n/m}\text{ is }o(t^{n/m})\text{ for }t\to\infty ,\\
\mu  _j(A_\gamma ^{-1})&-c_A^{m/n}j^{-m/n}\text{ is }o(j^{-m/n})\text{ for }j\to
\infty ,\endaligned
\end{equation*}
where 
\begin{equation*}
c_A=\int_{x\in \Omega ,\,a^0(x,\xi )<1}\,dx\GGdd\xi 
\end{equation*}
(and $c_0=c_A^{-m/n}$). The estimates have been shown for general $n$ and sharpened since then, with more
precision on the remainder, and the validity has been extended to
general elliptic operators $A$ and boundary conditions, and to
elliptic pseudodifferential operators $P$ on compact manifolds. These
improvements have a long history that we shall not try to account for
here in detail; they
are interesting not only because of the results but even more because
of the refined theories that were invented in connection with the
proofs (for example: Fourier integral operators). See e.g.\
H\"ormander \cite{GGH68}, \cite{GGH71}. Br\"uning \cite{GGB74}, Seeley \cite{GGSe78}, Ivrii
\cite{GGI82, GGI84, GGI91}, Safarov and Vassiliev \cite{GGSV97}.

We shall here be concerned with a slightly different question, namely
the spectral behavior of the difference between the resolvents of 
two realizations of $A$. 

Birman showed in \cite{GGB62} for second-order symmetric
uniformly strongly elliptic operators $A$ that the singular numbers $s_j(B)=\mu
_j(B^*B)^\frac12$ of the compact operator $B=\GGwA^{-1}-A_\gamma ^{-1}$ satisfy an upper estimate:
\begin{equation}
s_j(\GGwA^{-1}-A_\gamma ^{-1})\le Cj^{-2/(n-1)}, \text{ for all }j,\label{GG9.4}
\end{equation}
when $\GGwA $ is a selfadjoint
realization of $A$ defined by a Neumann or Robin  condition. In other
words, $\GGwA^{-1}-A_\gamma ^{-1}$ belongs $\mathfrak
S_{((n-1)/2)}$, where $\mathfrak S_{(p)}$ is  the space
of compact operators $B$ for which $s_j(B)$
is $O(j^{-\frac1p})$ (often called a weak Schatten class). It is particularly interesting
that Birman showed this not just for interior, but also for exterior
domains, and under low smoothness assumptions. 

We note in passing that the estimate \eqref{GG9.4} implies that
\begin{equation}
\GGwA^{-1}-A_\gamma ^{-1}\in {\mathcal C} _p\text{ for }p>(n-1)/2,\label{GG9.5}
\end{equation}
where  ${\mathcal C}_p$ is the  $p$-th Schatten class (consisting of the operators $B$ such that\linebreak
$(s_j(B))_{j\in{\mathbb N}}\in \ell_p({\mathbb N})$); but this is less
informative that \eqref{GG9.4}.

One of the fundamental ingredients in these studies is embedding
properties, more precisely the knowledge that an operator
$B$ that is continuous from $L^2(\Omega )$ to $H^s(\Omega )$ for some
$s>0$ ($\Omega $ bounded smooth $\subset {\mathbb R}^n$) is compact in
$L_2(\Omega )$ and
belongs to $\mathfrak
S_{(n/s)}$. However, this alone only gives upper 
bounds on the behavior of singular numbers. To get Weyl-type limit
properties one must know more about the differential or
pseudodifferential structure of the operators.

The estimate \eqref{GG9.4} was sharpened to a Weyl-type asymptotic
estimate in joint works of Birman and Solomiak 1978-80, see in particular \cite{GGBS80}
which showed a general principle for the spectrum of a ratio of two
quadratic forms, implying that 
\begin{equation}
s_j(\GGwA^{-1}-A_\gamma ^{-1})j^{2/(n-1)}\to c, \text{ for }j\to \infty ,\label{GG9.6}
\end{equation}
for interior and exterior smooth domains.

 Prior to this, a far-reaching result had been
 shown in \cite{GGG74}, Section 8:
We consider a symmetric, strongly elliptic $2m$-order operator acting in an
 $N$-dimensional vector bundle $E$ over a smooth
 compact Riemanninan manifold  $\GGcomega$ with boundary $\Sigma $, 
assuming that
 the Dirichlet realization $A_\gamma $ is invertible. Let $A_B$ be a
 selfadjoint invertible realization defined by a normal boundary condition
\begin{equation}
\sum_{k\le j}B_{jk}\gamma _ku=0,\; j=0,1,\dots, 2m-1,
\end{equation}
where the $B_{jk}$ are differential operators of order $j-k$ from
$E|_\Sigma $ to given vector bundles $F_j$ over $\Sigma $ (with $\dim
F_j\ge 0$); {\it normality}
of the boundary condition means that the $B_{jj}$ are {\it surjective
morphisms}. (That $A$ acts in the vectorbundle $E$ means that it is
locally $(N\times N)$-matrix-formed. In the scalar case, $N=1$ and the
$F_j$ are 0- or 1-dimensional, with $B_{jj}$ an invertible function
when $\dim F_j=1$. Ellipticity of the boundary condition requires in
particular that $\sum_j\dim F_j=mN$.) Denote $\oplus_{j>m}F_j=F^1$.

\begin{theorem}\label{GGTheorem 9.1} Let $T\colon V\to V$ be the operator
corresponding to $A_B$ by Theorem {\rm \ref{GGTheorem 5.1}}. There exists an {\bf isometry} $J\colon L_2(\Sigma
,F^1)\GGsimto V$ 
with inverse $J^{-1}=J^*$ (in the $\psi $dbo calculus), such that
\begin{equation*}
\mathcal T_1=J^*TJ
\end{equation*}
 acts like an elliptic invertible $\psi $do $\mathcal T$
 in $F^1$ of order $2m$, and $D(\mathcal T_1)=$ \linebreak$\{\varphi \in L_2(\Sigma
 ,F^1)\mid \mathcal T\varphi \in  L_2(\Sigma ,F^1)\}=H^{2m}(\Sigma
 ,F^1)$. Here $\mathcal T_1$ has the same spectrum as $T$, and its
 eigenvectors are mapped to the corresponding eigenvectors of $\mathcal
 T_1$ by the isometry $J$. Moreover,
\begin{equation}
A_B^{-1}-A_\gamma ^{-1}=\GGinj_VJ\mathcal T_1^{-1}J^*\GGpr_V,
\label{GG9.6a}\end{equation}
whereby the positive resp.\ negative eigenvalues satisfy
\begin{equation*}
\mu ^\pm_j(A_B^{-1}-A_\gamma ^{-1}) =\mu ^\pm _j({{\mathcal
    T}_1}^{-1})\text { for all }j.
\end{equation*}

It follows that with constants determined from the principal symbols,
\begin{equation}
\aligned
N^{\prime\pm}(t;A_B^{-1}-A_\gamma ^{-1})&=C^\pm t^{(n-1)/(2m)}
+O(t^{(n-1)/(2m)})\text{ for }t\to\infty ,\\
\mu ^\pm_j(A_B^{-1}-A_\gamma ^{-1})&=(C^\pm )^{2m/(n-1)}j^{-2m/(n-1)} +O(j^{-(2m+1)/(n-1)});
\endaligned
\label{GG9.7}
\end{equation}
here $N^{\prime\pm}(t;S)$ indicates the number of positive, resp.\ negative
eigenvalues of $S$ outside the interval $\,]-1/t,1/t[\,$. 

\end{theorem}

The two statements in \eqref{GG9.7} are equivalent, cf.\ e.g.\
\cite{GGG96}, Lemma A.5.
They follow from H\"o{}rmander \cite{GGH68}  for
elliptic $\psi $do's,
when the principal
symbol eigenvalues of $\mathcal T$ are simple; this restriction is
removed by results of Ivrii \cite{GGI82}.
See also Theorem \ref{GGTheorem 10.1} below.

We have recently checked that the proof extends to exterior domains,
for uniformly strongly elliptic systems.

The fine estimates with remainders depend on the ellipticity of the
$\psi $do.   For simple Weyl-type estimates,
Birman and Solomiak in
\cite{GGBS77} removed the ellipticity hypothesis, showing that 
\begin{equation*}
s_j(P)j^{k/n}\to c(p^0)\text{ for }j\to\infty 
\end{equation*} 
holds for any classical
$\psi $do $P$ of order $-k<0$ on a compact manifold of dimension
$n$. They even allowed a certain nonsmoothness of the homogeneous principal
symbol, both in $x$ and $\xi $, needing only a little more than
continuity. In the elliptic case, there are recent works of 
Ivrii dealing with remainder
estimates under weak smoothness hypotheses.
\medskip

In \cite{GGG84} we made an effort to increase the accessibility of the $\psi $dbo calculus by
publishing an introduction to it with several improvements, and showing
as a main result that any singular Green operator of negative order
and class 0 has a Weyl-type spectral asymptotics formula:

\begin{theorem}\label{GGTheorem 9.2}
Let $G$ be a  classical 
singular Green
operator of order $-k<0$ and class $0$  on an $n$-dimensional compact manifold
with boundary. It has a spectral asymptotics behavior
\begin{equation}
s_j(G)j^{k/(n-1)}\to c(g^0)\text{ for }j\to\infty ,\label{GG9.9}
\end{equation}
where $c(g^0)$ is a constant defined from on the principal symbol of
$G$. 
\end{theorem}

This was moreover used to show asymptotic formulas generalizing
\eqref{GG9.6}, both for interior and exterior domains, followed up in another
study \cite{GGG84a}  
including also the dependence on a spectral parameter $\lambda $.
Indeed, we have as an immediate corollary of Theorem \ref{GGTheorem 9.2}, 
also 
for nonselfadjoint cases:

\begin{corollary}\label{GGCorollary 9.3} Let $A$ be elliptic of order $2m$ with
invertible Dirichlet realization and let $A_B$ be an invertible realization
defined by a normal elliptic boundary condition. For any positive
integer $N$, $A_B^{-N}-A_\gamma
^{-N}$ is a singular Green operator of order $-2mN$ and class
$0$, and hence satisfies
\begin{equation}
s_j(A_B^{-N}-A_\gamma ^{-N})j^{2mN/(n-1) }\to c _N \text{
  for }j\to \infty ,\label{GG9.10}
\end{equation}
for a constant $c_N$ defined from the principal symbols.
\end{corollary}

Also exterior domains are considered in \cite{GGG84}, \cite{GGG84a}, where \eqref{GG9.10} is shown for realizations of
second-order operators and their iterates. The results apply of course
to resolvents by replacement of $A$ by $A-\lambda $; the $\lambda
$-dependence is studied in \cite{GGG84a}.
It is seen that the $\psi $dbo theory provides a forceful tool for such
questions, and we strongly recommend its use. 

\section{New spectral results} \label{GG10new}

Spectral estimates of resolvent differences have been taken up in
recent papers by  Behrndt et al.\ \cite{GGAB09, GGBLLLP10, GGBLL11} for second-order
operators and Malamud \cite{GGM10} for $2m$-order operators with 
normal boundary conditions, based on boundary triples methods. Here Schatten class and weak Schatten
class estimates are shown, relying on such estimates for
Sobolev space embeddings.

We have returned to the subject in \cite{GGG11} where we, besides 
showing new results
on perturbations of essential spectra, have reformulated and extended 
results in \cite{GGG84} on estimates like \eqref{GG9.10}, 
including general
differences and exterior domains. 
The central
ingredient is the estimate \eqref{GG9.9} for singular Green operators, plus
the fact that s.g.o.s give their essential contribution in a small
neighborhood of the boundary, also for exterior domains,
allowing cutoffs eliminating infinity.

An inspection of the results of \cite{GGG74} shows that the spectral
estimates in Theorem \ref{GGTheorem 9.1} can be further sharpened by
use of results of Ivrii \cite{GGI82}:

\begin{theorem}\label{GGTheorem 10.1} In the setting of 
  Theorem {\rm \ref{GGTheorem 9.1}}, assume in addition that the
  principal symbol of $\mathcal T$ satifies  Ivrii's
  conditions {\rm (H$_\pm$)} from {\rm \cite{GGI82}} (the
  bicharateristics through points of $T^*(\Sigma )\setminus 0$ are
  nonperiodic except for a set of measure zero). Then there are
  constants $C_1^\pm$ such that 
\begin{equation} 
N^{\prime\pm}(t;A_B^{-1}-A_\gamma ^{-1})=C^\pm t^{(n-1)/(2m)}+C_1^\pm t^{(n-2)/(2m)}
+o(t^{(n-2)/(2m)}).
\end{equation}
\end{theorem}

The proof is a direct application of \cite{GGI82} Th.\ 0.2 to
$\mathcal T$.

The formula \eqref{GG9.6a} is a special type of Kre\u\i{}n resolvent
formula with isometries, valid for selfadjoint realizations, but the  analysis in
\cite{GGG74} also implies Kre\u\i{}n formulas in the nonselfadjoint cases. Namely,
Th.\ 6.4 there shows how $T\colon V\to W$ is represented by a
realization $\mathcal L_1$
of a $\psi $do
$\mathcal L$ acting between vector bundles over $\Sigma
$, and here  $A_B$ is elliptic if and only if $\mathcal L$ is
elliptic (Cor.\ 6.10).  In the invertible elliptic case, formula \eqref{GG6.22} then
takes the form 
\begin{equation*}A_B^{-1}-A_\gamma
^{-1}=\GGinj_{V }\gamma _V^{-1}\Phi \mathcal L_1^{-1}\Psi ^*(\gamma
_W^*)^{-1}\GGpr_W,
\end{equation*}
with $\mathcal L_1$ acting like $-(B^{10}+B^{11}P_{\gamma ,\chi })\Phi
$ (notation explained in \cite{GGG74}); the right-hand side is
a composition of a Poisson operator, an elliptic $\psi
$do and the adjoint of a Poisson operator, all of mixed order. Its
$s$-numbers can be studied by reduction to an elliptic $\psi $do over
$\Sigma $, where Ivrii's sharp results
can be used.

Let us just demonstrate 
this for second-order operators,
for the formula \eqref{GG6.28}, with $L^{-1}$ elliptic of order $-1$:
Denote $\GGwA^{-1}-A_\gamma ^{-1}=S$. Then 
\begin{equation*}
\aligned
s_j(S)^2
&=s_j(K_{\gamma }L^{-1}(K'_{\gamma
  })^*)^2=
\mu _j(K_{\gamma }L^{-1}(K'_{\gamma
  })^*K'_{\gamma }(L^{-1})^*K_{\gamma
  }^*)\\
&=\mu _j(L^{-1}(K'_{\gamma })^*K'_{\gamma }(L^{-1})^*K_{\gamma
  }^*K_\gamma ),
\endaligned
\end{equation*}
where we used the general rule $\mu _j(B_1B_2)=\mu _j(B_2B_1)$.
Both operators  $P_1=K_{\gamma }^*K_{\gamma }$  and $P'_1=(K'_{\gamma
})^*K'_{\gamma }$
are selfadjoint positive elliptic $\psi $do's of order $-1$ (cf.\
e.g.\ \cite{GGG11b}, proof of Th.\ 4.4). Let $P_2=P_1^{\frac12}$, then we continue the
calculations as follows: 
$$
s_j(S)^2=\mu _j(L^{-1}P_1'(L^{-1})^*P_1)=\mu _j(P_2L^{-1}P_1'(L^{-1})^*P_2)=\mu _j(P_3),
$$
where $P_3=P_2L^{-1}P_1'(L^{-1})^*P_2$ is a selfadjoint positive  elliptic $\psi $do on
$\Sigma $ of order $-4$.
Applying Ivrii's theorem to $P_3^{-1}$, we conclude:

\begin{theorem}\label{GGTheorem 10.2} For the operator considered in
  Theorem {\rm \ref{GGTheorem 6.3}}, the $s$-numbers satisfy
\begin{equation} 
N'(t;\GGwA^{-1}-A_\gamma ^{-1})=C t^{(n-1)/2}+O(t^{(n-2)/2}).
\end{equation}
Moreover, if the
  principal symbol of $P_3^{-1}$ satifies  Ivrii's
  condition from {\rm \cite{GGI82}} (the
  bicharateristics through points of $T^*(\Sigma )\setminus 0$ are
  nonperiodic except for a set of measure zero), there is a
  constant $C_1$ such that 
\begin{equation} 
N'(t;\GGwA^{-1}-A_\gamma ^{-1})=C t^{(n-1)/2}+C_1 t^{(n-2)/2}
+o(t^{(n-2)/2}).
\end{equation}
\end{theorem}

Sharpened asymptotic formulas can also be obtained for differences between
re\-sol\-vents of two realizations that both differ from the Dirichlet
realization, by use of the analysis in \cite{GGG68} with a
general invertible realization $A_\beta $ as reference operator. 
\medskip

To give another example of applications of the $\psi $dbo theory, the following result is found
straightforwardly as a consequence of \cite{GGG84}:

\begin{theorem}\label{GGTheorem 10.3}
Let $A_{B }$ and $A_{B ' }$ be elliptic invertible
realizations of $A$ such that $B $ and $B '$ map into the same bundles
and have the {\bf same principal part}.
Then $A_{B }^{-1}-A_{B ' } ^{-1}$ is a singular Green
operator of order $-2m-1$ (since its principal part is zero), and
hence, by {\rm \eqref{GG9.9}},
\begin{equation*}
s _j(A_{B }^{-1}-A_{B ' } ^{-1})j^{(2m+1)/(n-1)}\to c \text{ for }j\to
\infty .
\end{equation*}
\end{theorem}

The singular Green operator will be of a still lower order $-2m-r$ if 
the first $r>1$ terms in the symbols of $B$ and $B'$ coincide.

\begin{example}\label{GGExample 10.4} 
As a  special case, we can compare two Robin conditions for a second-order operator $A$:
\begin{equation*}
\aligned &\GGwA_1 \text{ defined by }\nu _1u=b_1\gamma _0u,\\
&\GGwA_2 \text{ defined by }\nu _1u=b_2\gamma _0u;
\endaligned
\end{equation*}
$b_1, b_2\in C^\infty (\Sigma )$. When regarded from the point of view of
Theorem \ref{GGTheorem 10.3}, these are normal boundary
conditions $\nu _1u-B_i\gamma _0u=0$, where 
$B_1=b_1$ and $B_2=b_2$ considered as first-order operators
have principal part 0, so the boundary operators have the same
principal part. Then the s.g.o.\
$\GGwA_1^{-1}-\GGwA_2^{-1}$ is of order $-3$, and by \eqref{GG9.9},
\begin{equation}
s_j(\GGwA_1^{-1}-\GGwA_2^{-1})j^{3/(n-1)}\to c \text{ for }j\to\infty .\label{GG9.13} 
\end{equation}
\end{example}

\medskip
In \cite{GGBLLLP10},
 Berndt, Langer, Lobanov, Lotoreichik and Popov
prove upper estimates for this difference in the case $A=-\Delta -\lambda $, namely Schatten class estimates of
$s_j(\GGwA_1^{-1}-\GGwA_2^{-1})$ as in \eqref{GG9.5} with $(n-1)/2$ replaced
by $(n-1)/3$. 

In case the $b_i$ are $C^\infty $, the result is covered
by  \eqref{GG9.13} as explained above. However, the $b_i$ in
\cite{GGBLLLP10} are allowed to be  nonsmooth, namely to be in
$L_\infty(\Sigma)$, which goes outside the range covered by the smooth
$\psi $dbo theory. 

This led us to investigate how far we could push the proof of
asymptotic estimates \eqref{GG9.13} to make them valid for nonsmooth
choices of $b_i$. The outcome is published in \cite{GGG11b}, where it is
shown that \eqref{GG9.13} holds for symmetric second-order strongly
elliptic operators on smooth domains, when $b_1$ and $b_2$ are
piecewise $C^\varepsilon $ on $\Sigma $, having jumps at a smooth hypersurface.

The Schatten class estimates have been followed up by Berndt, Langer
and Lotoreichik in a study of selfadjoint realizations \cite{GGBLL11}.
\medskip

Unsolved questions of asymptotic estimates lie primarily in the
range of situations with limited smoothness. Resolvent formulas have
been studied in such general cases,  
\cite{GGPR09} and \cite{GGG08} for $C^{1,1}$-domains,
\cite{GGGM08, GGGM11} for Lipschitz and quasi-convex domains, 
\cite{GGAGW11} for a class of domains containing $C^{3/2+\varepsilon
}$, with a nonsmooth generalization of $\psi $dbo's. To our
knowledge, 
spectral asymptotic estimates have not yet been worked out for such
resolvent differences. Some upper
estimates are in selfadjoint cases known from Birman \cite{GGB62}.

A problem with a different flavor is the case of a {\it mixed
boundary condition}, such as prescribing for $-\Delta $ the 
Dirichlet condition on a
part $\Sigma _-$ of the boundary and a Neumann-type condition on the
other part $\Sigma _+$. Here there is a jump in the {\it order} of the
boundary condition. The domain of the realization is contained in
$H^{\frac32-\varepsilon }(\Omega )$ only for $\varepsilon >0$, so it
is not covered by boundary triples methods requiring the domain to be
in $H^{\frac32}(\Omega )$. Spectral upper estimates are known from
\cite{GGB62}. A spectral asymptotic estimate was
obtained recently in \cite{GGG11c}, based on somewhat intricate
applications of results on nonstandard pseudodifferential operators.

There are also other questions that can benefit from
pseudodifferential methods, for example the
study of spectral asymptotics of the nonelliptic 
Kre\u\i{}n-like extensions, cf.\
\cite{GGG11a}. 




\begin{thebibliography}{100}

\bibitem[AGW11]{GGAGW11}
 H. Abels, G. Grubb and I. Wood.  Extension theory and  Kre\u\i{}n-type resolvent
  formulas for nonsmooth boundary value
  problems,  arXiv:1008.3281, to appear.


\bibitem[A65]{GGA65}
S. Agmon. {\em Lectures on Elliptic Boundary  Value
Problems}. Van Nostrand Math.\ Studies, D.\ Van Nostrand Publ.\
Co., Princeton 1965.
 


\bibitem[ADN59]{GGADN59}
S. Agmon, A. Douglis and L. Nirenberg. 1964.
Estimates near the boundary for solutions of elliptic partial
differential equations satisfying general boundary conditions, I.
{\em  Comm.\ Pure Appl.\
Math}, {\bf l2}, 623--727.
 

\bibitem[AB09] {GGAB09} D. Alpay and J. Behrndt. 2009. Generalized
  Q-functions and Dirichlet-to-Neumann maps for elliptic differential
  operators. {\em J. Funct. Anal.} {\bf 257}, 1666--1694.


\bibitem[AP04]   {GGAP04}   
W. O. Amrein and D. B. Pearson. 2004.    $M$ operators: a
generalisation of Weyl-Titchmarsh theory.  {\em   
J.~Comp.~Appl. Math.} {\bf    171},  1--26.  

\bibitem[A99]  {GGA99}  
Yu. M. Arlinskii. 1999.  On functions connected with
sectorial operators and their extensions. {\em  Integral Equations
Operator Theory} {\bf  33},  125--152. 
 

 

 

 \bibitem[BL07]   {GGBL07}   
J.~Behrndt and M.~Langer.  2007.     Boundary value
problems for elliptic partial  differential operators on bounded
domains. {\em    J.~Funct.~Anal.}  {\bf   243},  536--565.

\bibitem[BLLLP10]  {GGBLLLP10}
J. Behrndt, M. Langer, I. Lobanov,
V. Lotoreichik and I. Popov. 2010.    A remark on Schatten-von Neumann
properties of resolvent differences of generalized Robin Laplacians on
bounded domains.
{\em  J. Math. Anal. Appl.}
{\bf  371},   750-758. 

 
\bibitem[BLL11]  {GGBLL11}
J. Behrndt, M. Langer and Lotoreichik. Spectral estimates for
differences of resolvents of selfadjoint elliptic operators.
 arXiv:1012.4596, to appear.


\bibitem[B56]   {GGB56}   
M. S. Birman. 1956.   On the theory of self-adjoint
extensions of positive definite operators. {\em  
Mat. Sb. N.S.} {\bf  38(80)},   431--450. (In Russian.)  

\bibitem[B62]  {GGB62}
M. S. Birman. 1962. Perturbations of the continuous
spectrum of a singular elliptic operator by varying the boundary and
the boundary conditions.
{\em  Vestnik Leningrad. Univ.} {\bf  17},   22--55.
 English translation in
{\em Spectral theory of differential operators}, 
 Amer. Math. Soc. Transl. Ser. 2, {\bf 225}. Amer. Math. Soc.
Providence, R.I.  2008,   19--53.  
  
 

\bibitem[BS77]  {GGBS77}  
 M. S. Birman and M. Z.  Solomyak. 1977.  Asymptotic
behavior of the spectrum of pseudodifferential operators with
anisotropically homogeneous symbols. {\em Vestnik
Leningrad. Univ.} {\bf 13}, 13--21.  English
translation in {\em Vestn. Leningr. Univ. Math.} {\bf 10} (1982),
 237--247. 


\bibitem[BS80]     {GGBS80}     
M. S. Birman and M. Z. Solomyak. 1980.
  Asymptotics of the spectrum of variational problems on
solutions of elliptic equations in unbounded domains. {\em 
  Funkts. Analiz Prilozhen.}
  {\bf 14},   27--35. English translation in
{\em  Funct. Anal. Appl.} {\bf 14} (1981),   267--274.
 
 

\bibitem[B71]  {GGB71}  
  L.~Boutet de Monvel. 1971.    Boundary problems for pseudodifferential
operators, {\em   
 Acta Math.} {\bf 126},   11--51. 
 
\bibitem[BGW09]  {GGBGW09}  
B. M. Brown, G. Grubb, and I. G. Wood. 2009.   $M$-functions for closed
extensions of adjoint pairs of operators with applications to elliptic boundary
problems. {\em  Math. Nachr. } {\bf  282},  314--347.
   

\bibitem[BMNW08]  {GGBMNW08}  
B.~M.~Brown, M.~Marletta,  S.~Naboko and
 I.~G.~Wood. 2008.    
Boundary triplets and M-functions for non-selfadjoint operators, with
applications to elliptic PDEs and block operator matrices. {\em   
J. Lond. Math. Soc.} {\bf  77}, 700--718.


\bibitem[B74]  {GGB74}
J. Br\"u{}ning. 1974.   Zur Absch\"a{}tzung der
Spektralfunktion elliptischer Operatoren. {\em  Math. Z.} {\bf 137},
  75--85.  

\bibitem[BGP06]   {GGBGP06}   
J.~Br\"{u}ning, V.~Geyler and K.~Pankrashkin. 2008.  
Spectra of self-adjoint extensions and applications to solvable
Schr\"{o}dinger operators. {\em  Rev. Math. Phys.} {\bf   20}, 
1--70.  

\bibitem[CZ57]  {GGCZ57}
A.-P. Calder\'o{}n and A. Zygmund. 1957.   Singular integral operators and
differential equations. {\em  Amer. J. Math.} {\bf  79},   901--921.
 

\bibitem[CH53]  {GGCH53}
R. Courant and D. Hilbert. {\em Methods of
Mathematical Physics, Vol. I}. Interscience Publishers, Inc., New York, N.Y. 1953.
 

\bibitem[CH62]  {GGCH62}  
R. Courant and D. Hilbert. {\em Methods of
Mathematical Physics, Vol. II: Partial Differential Equations} (by
R. Courant). Interscience
Publishers (a division of John Wiley \& Sons), New
York-London 1962.
  

\bibitem[DM91]   {GGDM91}   
V. A.~Derkach and M. M.~Malamud. 1991.   Generalized resolvents and
the boundary value problems for Hermitian operators with gaps. {\em   J.~Funct.~Anal.}  {\bf   95},   1--95. 

\bibitem[F34] {GGF34}
  K.\ Friedrichs. 1934.
  Spektraltheorie halbbeschr\"ankter Operatoren und Anwendung
auf die Spektralzerlegung von Differentialoperatoren.
{\em  Math. Ann.}
{\bf  109},  465--487.
 

\bibitem[GM08]   {GGGM08} 
F. Gesztesy and M. Mitrea. 2008.   Generalized Robin boundary
conditions, Robin-to-Dirichlet maps, and Krein-type resolvent formulas
for Schr\"odinger operators on bounded Lipschitz domains. {\em  Perspectives in Partial Differential Equations, Harmonic Analysis
and Applications: A Volume in Honor of Vladimir G. Maz'ya's 70th Birthday,}
Proceedings of
Symposia in Pure Mathematics (eds. D. Mitrea and M. Mitrea)  
 {\bf  79}, Amer. Math. Soc., Providence,
R.I.,  105--173.
 
\bibitem[GM11]   {GGGM11} 
F.~Gesztesy and M.~Mitrea. 2011.
A description of all selfadjoint extensions of the Laplacian
 and Kre\u\i{}n-type
resolvent formulas in nonsmooth domains. {\em  J. Analyse Math.}
{\bf 113}, 53--172.

\bibitem[GM76]  {GGGM76}
M. L. Gorbachuk and V. A. Mikhailets. 1976.  Semibounded
selfadjoint extensions of symmetric operators. {\em  Dokl. Akad. Nauk
SSSR} {\bf  226}. English translation in {\em 
Soviet Math. Doklady} {\bf  17},  185--186.
 

 \bibitem[GG91]   {GGGG91}
V. I.~Gorbachuk and M. L.~Gorbachuk.
{\em Boundary value problems for operator differential equations}.
Kluwer, Dordrecht 1991.
  

\bibitem[G68]  {GGG68}
G. Grubb. 1968.
  A characterization of the non-local boundary value problems
associated with an elliptic operator
{\em  Ann.\ Scuola Norm.\ Sup.\ Pisa}
{\bf 22}, 425--513.
 

\bibitem[G70]  {GGG70}
G.~Grubb. 1970.   Les probl\`emes aux limites g\'en\'eraux d'un
op\'erateur elliptique, provenant de la th\'eorie variationnelle.
 {\em Bull.~ Sc.~Math.} {\bf 94},  113--157. 

\bibitem[G71]  {GGG71}
G.~Grubb. 1971.   On coerciveness and semiboundedness of general boundary
value problems. {\em  
Israel J. Math.} {\bf 10},  32--95. 

\bibitem[G73] {GGG73}
  G. Grubb. 1973. 
  Weakly semibounded boundary problems and
sesquilinear forms.
{\em  Ann. Inst. Fourier Grenoble}
{\bf  23}, 145--194.
 

\bibitem[G74]  {GGG74}
G. Grubb. 1974.  Properties of normal boundary problems for elliptic
even-order systems. {\em  Ann.\ Scuola Norm.\ Sup.\ Pisa} {\bf 1}
(ser.IV), 1--61.
 


\bibitem[G83] {GGG83}
  G. Grubb. 1983.
  Spectral asymptotics for the ``soft'' selfadjoint
extension of a symmetric elliptic differential operator.
{\em  J. Operator
Theory}
{\bf 10}, 9--20.
 


\bibitem[G84] {GGG84}
  G. Grubb. 1984.
  Singular Green operators and their spectral asymptotics.
{\em  Duke Math. J.}
{\bf  51},  477--528.
 

\bibitem[G84a]  {GGG84a}
G. Grubb. 1984.   Remarks on trace extensions for
exterior boundary problems
{\em  Comm. Partial Diff. Equ.}  {\bf  9},   231--270.
 


 \bibitem[G96]  {GGG96}
G.~Grubb. {\em Functional Calculus of Pseudodifferential
     Boundary Problems},
 Pro\-gress in Math.\ vol.\ 65, Second Edition.  Birkh\"auser,
  Boston  1996. 

\bibitem[G08]  {GGG08}
G. Grubb. 2008. Krein resolvent formulas for elliptic boundary
problems in nonsmooth domains. {\em Rend. Sem. Mat. Univ. Pol. Torino}
{\bf 66}, 13--39. 

\bibitem[G09]  {GGG09}
G. Grubb. {\em  Distributions and operators.} Graduate
Texts in Mathematics {\bf 252}. Springer, New York 2009.
  
 
\bibitem[G11]    {GGG11}
G. Grubb. 2011. 
Perturbation of essential spectra of exterior elliptic problems.
{\em  Applicable Analysis} {\bf 90}, 103--123.
  


\bibitem[G11a]  {GGG11a}
G. Grubb.   Krein-like extensions and the lower boundedness problem for elliptic operators on exterior domains. 
arXiv:1002.4549, to appear.
 

\bibitem[G11b]  {GGG11b}
G. Grubb. 2011.   Spectral asymptotics for Robin
problems with a discontinuous coefficient. {\em J. Spectral Theory}
{\bf 1}, 155--177.

\bibitem[G11c] {GGG11c}
G. Grubb. 2011.
The mixed boundary value problem, Krein resolvent formulas and
spectral asymptotic estimates. {\em   J. Math. Anal. Appl.} {\bf 382},
339--263.  
 
 

\bibitem[H63]  {GGH63}
L.\ H\"ormander. {\em Linear Partial Differential Operators,}
Grundlehren Math. Wiss. vol. 116. Springer Verlag, Berlin 1963.
 

\bibitem[H65]  {GGH65}
L. H\"o{}rmander. 1965.  Pseudo-differential operators. {\em  Comm.\ Pure
Appl.\ Math.} {\bf 18},  501--517.
 

\bibitem[H68]  {GGH68}
L. H\"o{}rmander. 1968.  The spectral function of an elliptic
operator. {\em  Acta Math.} {\bf 121}, 193--218.
 

\bibitem[H71]  {GGH71}
L. H\"o{}rmander. 1971.  Fourier integral operators I. {\em  Acta
Math.} {\bf 127}, 79--183.
 
\bibitem[H85] {GGH85}
  L. H\"ormander. {\em  The Analysis of Linear Partial Differential Operators III,
Pseu\-do-dif\-fe\-ren\-tial Operators,} Grundlehren Math. Wiss. vol. 274. 
Sprin\-ger Ver\-lag, Berlin
  1985.
 

\bibitem[I82]  {GGI82}  
V. Ja. Ivrii. 1982. Accurate spectral
  asymptotics for elliptic operators that act in vedtor bundles.
{\em Functional Analysis Prilozhen} {\bf 16} no. 2, 30--38, English
translation in {\em Functional Analysis Appl.} {\bf 16} (1983), 101--108.

\bibitem[I84]  {GGI84}  
V. Ivrii. {\em Precise spectral asymptotics for
elliptic 
operators acting in fiberings over manifolds with boundary}.
Lecture Notes in Mathematics, 1100. Springer-Verlag, Berlin
 1984.

\bibitem[I91]  {GGI91}
V. Ivrii. {\em Microlocal Analysis and Precise Spectral Asymptotics}.
Springer-Verlag, Berlin
 1991.  

\bibitem[KN65]   {GGKN65} 
J. J. Kohn and L. Nirenberg. 1965.   An algebra of
pseudo-differential operators. {\em  Comm. Pure Appl. Math.} {\bf 18},
 269--305. 
 

\bibitem[K75]   {GGK75}
A. N.~Ko\v{c}ube\u{\i}. 1975.  {Extensions of symmetric operators and
  symmetric binary relations}. {\em   Math.~Notes} (1) {\bf   17},   25--28.

\bibitem[KK04]   {GGKK04}
N.D.~Kopachevski{\u{\i}} and S.G.~Kre\u{\i}n. 2004.
   Abstract Green formula for a triple of Hilbert spaces, abstract boundary value and spectral problems.
 {\em  Ukr.~Math.~Bull.} {\bf  1},   77--105.
 

\bibitem[K47]   {GGK47}
M. G. Kre\u\i{}n. 1947.   The theory of self-adjoint extensions of
semi-bounded Hermitian transformations and its applications. I.
{\em  Mat.\ Sbornik } {\bf  20}  431--495. (In Russian.) 
 

\bibitem[L85]  {GGL85}
O. Ladyshenskaya. {\em The Boundary Value Problems of
Mathematical Physics}. Springer-Verlag, New York 1985.
 

\bibitem[LM68]   {GGLM68}
J.-L. Lions and E. Magenes. {\em  Probl\`emes aux
limites non homog\`enes et applications, 1}. 
 \'Editions Dunod, Paris 1968.
 

\bibitem[L53]  {GGL53}
   Ya. B. Lopatinski\u\i{}. 1953.   On a method of
reducing boundary problems for a system of differential equations of
elliptic type to regular integral equations. {\em  Ukrain. Mat. Zb.} {\bf 
5},  123--151. 
 

 \bibitem[LS83]  {GGLS83}
V. E.~Lyantze and  O. G.~Storozh. {\em Methods of the Theory
of Unbounded
Operators}. Naukova Dumka,  Kiev 1983. 
 (In   Russian.) 
 

\bibitem[M10]   {GGM10}
M. M. Malamud.   Spectral theory of elliptic
operators in exterior domains. {\em  Russian J. Math. Phys.} {\bf  17},
  96--125.  

\bibitem[MM02]  {GGMM02}
{M. M. Malamud} and  {V. I. Mogilevskii.} 2002.  Kre\u\i n
  type formula for canonical resolvents of dual pairs of linear
  relations. {\em  
Methods Funct. Anal. Topology} (4) {\bf 8},  72--100.

\bibitem[M48]   {GGM48}   
S. G. Mihlin. 1948.   Singular integral equations.
{\em  Uspehi
Matem. Nauk (N.S.)} {\bf  3},  29--112. English translation
in Amer. Math. Soc. Translations {\bf 24} (1950), 116 pp.
 
 
\bibitem[N29]   {GGN29}   
J. von Neumann. 1929.   Allgemeine Eigenwerttheorie Hermitescher 
Funktionaloperatoren. {\em  Math. Ann.} {\bf  102},   49--131.
 
\bibitem[N55] {GGN55}
  L. Nirenberg. 1955.   Remarks on strongly elliptic partial
differential equations. {\em  Comm. Pure Appl. Math.} {\bf  8},
  649--675.
 
\bibitem[PR09]  {GGPR09}  
A. Posilicano and L. Raimondi. 2009.  Krein's resolvent
formula for self-adjoint extensions of symmetric second-order elliptic
differential operators. {\em J. Phys. A} {\bf 42}  015204,
11 pp. 


\bibitem[RS82]  {GGRS82}  
S.\ Rempel and B.-W.\ Schulze. {\em Index Theory of
Elliptic Boundary Problems}. A\-ka\-de\-mie-Verlag,
Berlin 1982.
 

\bibitem[R07]   {GGR07}   
V. Ryzhov. 2007.  A general boundary value problem and
its Weyl function.  {\em  Opuscula Math.} {\bf  27}   305--331. 
 


\bibitem[SV97]  {GGSV97}  
Yu. Safarov and D. Vassiliev. {\em The Asymptotic
Distribution of Eigenvalues of Partial Differential
Operators.} Translated from the Russian manuscript by the
authors. Translations of Mathematical Monographs, {\bf 155}. American
Mathematical Society, Providence, R.I.  1997.
  

\bibitem [Sch50] {GGSch50}
  L. Schwartz. {\em  Th\'eorie des distributions I--II}. 
Hermann, Paris
 1950--51.
 


\bibitem[Se65]   {GGSe65}   
R. T. Seeley. 1965.  Refinement of the functional calculus of
Calderon and Zygmund. {\em  Proc. Konikl. Nederl. Akad. Wetensch.} {\bf 
68},   521--531.
 

\bibitem[Se78]   {GGSe78}   
R. T.  Seeley. 1978.  A sharp asymptotic remainder
estimate for the eigenvalues of the Laplacian in a domain of $R^{3}$.
{\em  Adv. in Math.} {\bf  29},  244--269.
 

\bibitem[Sh53]   {GGSh53}   
Z. Ya. Shapiro. 1953.   On general boundary value
problems of elliptic type. {\em  Isz. Akad. Nauk, Math. Ser.}
{\bf  17},   539--562. 

\bibitem[So50]  {GGSo50}  
S. L.  Sobolev. {\em  Some applications of functional
analysis in mathematical physics}.
Izdat. Leningrad. Gos. Univ., Leningrad  1950. English translation
by F. E. Browder. Translations
of Mathematical Monographs, {\bf 7}, American Mathematical Society,
Providence, R.I. 1963.
 

\bibitem[T80]  {GGT80}  
F. Treves. {\em Introduction to Pseudodifferential and Fourier
Integral Operators, 1-2.} Plenum Press, New
York 1980.
 

 \bibitem[V80]   {GGV80}   
{L. I.~Vainerman}. 1980.    On extensions of closed operators in
Hilbert space.
{\em Math.~Notes} {\bf  28},   871--875. 

\bibitem[V52] {GGV52} 
M. I.~Vishik. 1952.   {On general boundary value problems for elliptic
differential operators}. {\em  
Trudy Mosc. Mat. Obsv.}  {\bf 1},  187--246.
English translation in {\em  Amer. Math. Soc. Transl. (2) } {\bf 24}
  (1963) 107--172. 

\bibitem[W12]  {GGW12}  
H. Weyl. 1912.   Das asymptotische Verteilungsgesetz der
Eigenwerte linearer partieller Differentialgleichungen (mit einer
Anwendung auf die Theorie der Hohlraumstrahlung). {\em 
Math. Ann.} {\bf 71},   441-479.
 

\end{thebibliography}
\end{document}